\documentclass{gtart_h}  

\def\ifplaintex{\expandafter\ifx\csname documentclass\endcsname\relax}

\ifplaintex 
\else
\fi

\expandafter\ifx\csname beginpicture\endcsname\relax
\expandafter\ifx\csname documentclass\endcsname\relax
\input pictex \else
\input prepictex \input pictex \input postpictex \fi\fi

\def\gt{{\mathsurround=0pt\it $\cal G\mskip-2mu$eometry \&\ 
$\cal T\!\!$opology}}        %  journal title in recommended style

\def\gtp{{\mathsurround=0pt\it $\cal G\mskip-2mu$eometry \&\ 
$\cal T\!\!$opology $\cal P\!$ublications}}  % GT publications

\def\lognumber#1{\def\thelognumber{#1}}
\def\volumenumber#1{\def\thevolumenumber{#1}}
\def\papernumber#1{\def\thepapernumber{#1}}
\def\volumeyear#1{\def\thevolumeyear{#1}}

\def\pagenumbers#1#2{\def\startpage{#1}\def\finishpage{#2}}
\def\published#1{\def\publishdate{#1}}
\def\proposed#1{\def\theproposer{#1}}
\def\seconded#1{\def\theseconders{#1}}
\def\received#1{\def\receiveddate{#1}}
\def\revised#1{\def\reviseddate{#1}}
\def\accepted#1{\def\accepteddate{#1}}

\long\def\asciiabstract#1{\long\def\theasciiabstract{#1}}

\let\\\par\let\thelognumber\relax
\let\thevolumenumber\relax\let\thepapernumber\relax
\let\thevolumeyear\relax\let\thesamplenumber\relax\let\startpage\relax
\let\finishpage\relax\let\publishdate\relax\let\receiveddate\relax
\let\reviseddate\relax\let\accepteddate\relax\let\theasciititle\relax
\let\theasciiauthors\relax
\let\theasciiabstract\relax
\let\theasciiemail\relax\let\theshortauthors\relax\let\theshorttitle\relax

\long\def\maketitlep{   % start of definition of \maketitlep

\count0=\startpage

\gt\hfill      %   Journal title (top left) 
\beginpicture
\setcoordinatesystem units <0.33truein, 0.33truein> point at 2.2 0.9
\setplotsymbol ({$\cal G$})
\plotsymbolspacing=9truept
\circulararc 315 degrees from 0 1 center at 0 0
\setplotsymbol ({$\cal T$})
\circulararc 315 degrees from 1 -1 center at 1 0
\endpicture
\break
{\small\ifx\thesamplenumber\relax % sample?  
Volume \else Sample
\fi\thevolumenumber\ (\thevolumeyear)
\startpage--\finishpage\nl
Published: \publishdate}
\vglue 0.5truein plus 0.4fil minus 0.1truein

{\parskip=0pt\leftskip 0pt plus 1fil\def\\{\par\smallskip}{\ifplaintex\large
\else\Large\fi\bf\thetitle}\par\medskip}   

\vglue 0pt plus 0.1fil 

{\parskip=0pt\leftskip 0pt plus 1fil\def\\{\par}{\sc\theauthors}
\par\medskip}

\vglue 0pt plus 0.1fil 

{\small\parskip=0pt\let\newline\\
{\leftskip 0pt plus 1fil\def\\{\par}{\sl\theaddress}\par}
\expandafter\ifx\theemail\relax    % email address?
\relax\else\vglue 5pt plus 0.02fil minus 2pt\def\\{\stdspace{\rm 
and}\stdspace} 
\cl{Email:\stdspace\tt\theemail}\fi
\ifx\theurl\relax                  % URL given?
\relax\else\vglue 5pt plus 0.02fil minus 2pt\def\\{\stdspace{\rm 
and}\stdspace}
\cl{URL:\stdspace\tt\theurl}\fi\par}

\vglue 7pt plus 0.3fil minus 3pt

{\bf Abstract}
\vglue 5pt plus 0.1fil minus 2pt

\theabstract

\vglue 7pt plus 0.3fil minus 3pt

{\bf AMS Classification numbers}\quad Primary:\quad \theprimaryclass

Secondary:\quad \thesecondaryclass

\vglue 5pt plus 0.3fil minus 2pt

{\bf Keywords:}\quad \thekeywords

\vglue 10pt plus 0.5fil minus 5pt

{\small  Proposed: \theproposer\hfill Received: \receiveddate\nl
Seconded: \theseconders\hfill 
\ifx\reviseddate\relax                         % paper revised?
Accepted: \accepteddate                        % no
\else
Revised: \reviseddate                          % yes
\fi}
\eject
}       %  end of definition of \maketitlep

\let\maketitlepage\maketitlep
\let\maketitle\maketitlepage

\font\phead=cmsl9 scaled 950
\font=cmsl9 scaled 1050
\font\pnum=cmbx10 scaled 913
\font\lnum=cmbx10 
\font\pfoot=cmsl9 scaled 950
\font=cmsl9 scaled 1050
\ifplaintex
\headline{\vbox to 0pt{\vskip -4.5mm\line{\small\phead\ifnum
\count0=\startpage ISSN 1364-0380 (on line)
1465-3060 (printed) \hfill {\pnum\folio}\else\ifodd\count0\def\\{ }% 
\ifx\theshorttitle\relax\thetitle\else\theshorttitle\fi\hfill{\pnum\folio}
\else\def\\{ and }{\pnum\folio}\hfill\ifx\theshortauthors\relax\theauthors
\else\theshortauthors\fi\fi\fi}\vss}}
\footline{\vbox to 0pt{\vglue 0mm\line{\small\pfoot\ifnum\count0=\startpage
\copyright\ \gtp\hfill\else
\gt, Volume \thevolumenumber\ (\thevolumeyear)\hfill\fi}\vss
}}
\else
\makeatletter
\def\@oddhead{{\small\lhead\ifnum\count0=\startpage ISSN 1364-0380 (on line)
1465-3060 (printed) \hfill {\lnum\number\count0}\else\ifodd\count0
\def\\{ }\ifx\theshorttitle\relax \thetitle \else\theshorttitle\fi\hfill
{\lnum\number\count0}\else\def\\{ and }{\lnum\number\count0}
\hfill\ifx\theshortauthors\relax 
\theauthors\else\theshortauthors\fi\fi\fi}}\def\@evenhead{\@oddhead}
\def\@oddfoot{\small\lfoot\ifnum\count0=\startpage\copyright\ \gtp\hfill\else
\gt, Volume \thevolumenumber\ (\thevolumeyear)\hfill\fi}
\def\@evenfoot{\@oddfoot}
\makeatother
\fi

\newwrite\gtoutfile
\long\gdef\makeheadfile{  %%% start of definition of \makeheadfile
{\def\\{, }\def\s{ }
\immediate\openout\gtoutfile head.xxx
\immediate\write\gtoutfile{Proxy-for: \ifx\theasciiauthors\relax
\theauthors\else\theasciiauthors\fi\s<\ifx\theasciiemail\relax\theemail\else\theasciiemail\fi>}
\immediate\write\gtoutfile{\noexpand\\}
\immediate\write\gtoutfile{Authors: \ifx\theasciiauthors\relax
\theauthors\else\theasciiauthors\fi}
{\def\\{ }\immediate\write\gtoutfile{Title: \ifx\theasciititle\relax
\thetitle\else\theasciititle\fi}}
\immediate\write\gtoutfile{Subj-class: GT or SG or MG etc}
\immediate\write\gtoutfile{MSC-class: \theprimaryclass\ifx\thesecondaryclass\relax\else, \thesecondaryclass\fi}
\immediate\write\gtoutfile{Journal-ref: Geom. Topol. \thevolumenumber
(\thevolumeyear) \startpage-\finishpage}
\immediate\write\gtoutfile{Comments: Published by Geometry and Topology at}
\immediate\write\gtoutfile{\s\s http://www.maths.warwick.ac.uk/gt/GTVol\thevolumenumber/paper\thepapernumber.abs.html}
\immediate\write\gtoutfile{\noexpand\\}
\immediate\write\gtoutfile{}
\ifx\theasciiabstract\relax
\immediate\write\gtoutfile{\theabstract}\else
\immediate\write\gtoutfile{\theasciiabstract}\fi
\immediate\write\gtoutfile{}
\immediate\write\gtoutfile{\noexpand\\}
\immediate\write\gtoutfile{}
\immediate\closeout\gtoutfile}}  %%% end of definition of \makeheadfile

\def\maketitlepage{\maketitlep\makeheadfile}
\let\maketitle\maketitlepage

\lognumber{472}

\received{25 June 2004}
\volumenumber{9}\papernumber{21}\volumeyear{2005}
\pagenumbers{935}{970}   
\revised{24 September 2004}
\published{25 May 2005}
\accepted{18 January 2005}
\proposed{Leonid Polterovich}
\seconded{Yasha Eliashberg, Tomasz Mrowka}

\usepackage{amsmath,graphicx,amssymb,amscd}

\newtheorem{theorem}{Theorem}[section]
\newtheorem{proposition}[theorem]{Proposition}
\newtheorem{corollary}[theorem]{Corollary}
\newtheorem{lemma}[theorem]{Lemma}
\newtheorem*{thm*}{Theorem}
\newtheorem*{prop*}{Proposition}
\newtheorem*{lem*}{Lemma}
\theoremstyle{remark}
\newtheorem{definition}[theorem]{Definition}
\newtheorem{remark}[theorem]{Remark}

\newcommand{\tmop}[1]{\ensuremath{\operatorname{#1}}}
\newcommand{\tmem}[1]{{\em #1\/}}
\newcommand{\tmdfn}[1]{\emph{#1}}

\newcommand{\sft}{\mathcal{S}^\text{soft}}
\newcommand{\sftio}{\left( {\sft} \right)_{\infty,0}}
\newcommand{\sio}{\mathcal{S}_{\infty,0}}

\newcommand{\s}{S^2}

\newcommand{\acs}{almost complex structure}

\begin{document}

\title{Symplectomorphism groups and isotropic skeletons}
\authors{Joseph Coffey}
\address{Courant Institute for the Mathematical Sciences, New York
University\\251 Mercer Street, New York, NY 10012, USA}
\email{coffey@cims.nyu.edu}

\begin{abstract}
  The symplectomorphism group of a 2--dimensional surface is homotopy
  equivalent to the orbit of a filling system of curves. We give a
  generalization of this statement to dimension 4. The filling system of
  curves is replaced by a decomposition of the symplectic 4--manifold $( M,
  \omega )$ into a disjoint union of an isotropic 2--complex $L$ and a disc
  bundle over a symplectic surface $\Sigma$ which is Poincare dual to a
  multiple of the form $\omega$. We show that then one can recover the
  homotopy type of the symplectomorphism group of $M$ from the orbit of the
  pair $( L, \Sigma ) .$ This allows us to compute the homotopy type of
  certain spaces of Lagrangian submanifolds, for example the space of
  Lagrangian $\mathbb{RP}^2 \subset \mathbb{CP}^2$ isotopic to the standard
  one.
\end{abstract}
\asciiabstract{%
The symplectomorphism group of a 2-dimensional surface is homotopy
equivalent to the orbit of a filling system of curves. We give a
generalization of this statement to dimension 4.  The filling system
of curves is replaced by a decomposition of the symplectic 4-manifold
(M, omega) into a disjoint union of an isotropic 2-complex L and a
disc bundle over a symplectic surface Sigma which is Poincare dual to
a multiple of the form omega.  We show that then one can recover the
homotopy type of the symplectomorphism group of M from the orbit of
the pair (L, Sigma).  This allows us to compute the homotopy type of
certain spaces of Lagrangian submanifolds, for example the space of
Lagrangian RP^2 in CP^2 isotopic to the standard one.}

\maketitlepage

\section{Introduction}

\subsection{Surface diffeomorphisms}

A finite set $L = \{ \gamma_i \}$ of simple, closed, transversally
intersecting parametrized curves on a surface $X$ \emph{fills} if $X
\backslash \{ \gamma_i \}$ consists solely of discs. Endow $X$ with a
symplectic structure $\omega$, and denote by $\mathcal{L}$ the orbit of $L$
under the action of the symplectomorphism group $\mathcal{S}$. We have the
following classical result:

\begin{theorem}
  \label{thm:classical sym}Let $( X, \omega )$ be a symplectic surface. Then
  the orbit map $\mathcal{S} \to \mathcal{L}$ is a homotopy equivalence.
\end{theorem}

\begin{proof}[Sketch proof]
  We must show that the stabilizer of $\{ \gamma_i \}$ in
  $\mathcal{S}$ (the symplectomorphisms of a disjoint union of
  discs, fixing their boundaries) is contractible.
    Moser's lemma allows us to replace symplectomorphisms with
    diffeomorphisms. (The inclusion is a weak deformation retract.)
    By the Riemann mapping theorem, these act transitively on the space
    of complex structures on the disc which are standard at the boundary.
    This set of complex structures is contractible, and thus the
    homotopy type is reduced to that of the complex automorphisms of the disc
    which fix the boundary, a contractible set.
\end{proof}

In this paper we will prove a 4 dimensional analog of this statement. The
proof follows a similar outline, although it is of course more difficult to
carry out each step. Our tool is again the Riemann mapping theorem, but this
time augmented by the theory of $J$--holomorphic spheres in sphere bundles over
surfaces, developed by Gromov, Lalonde and McDuff
{\cite{phol,Lalonde,Sympalmostcomp}}.

\subsection{Biran decompositions: a higher dimensional analog of filling
systems of curves}

Paul Biran {\cite{biran:Lspine}} recently showed that every K\"ahler manifold
$M$ whose symplectic form lies in a rational cohomology class admits a
decomposition
\[ {\textstyle M = L \coprod E} \]
where $L$ is an embedded, isotropic cell complex and $E$ is a symplectic disc
bundle over a hypersurface $\Sigma$. $L$ is called an isotropic skeleton of
$M$.

We will argue that a Biran decomposition of a symplectic $4$--manifold should
be regarded as the $4$--dimensional analog of a filling system of curves.
Indeed, when $M$ is a surface $L$ is a filling system of curves, $\Sigma$ is a
union of points (one in each disc inside $M \backslash L$) and $E$ is the
union of discs.

In higher dimensions we have less understanding of the possible singularities
of the spine $L$ and as a result we prove a weaker, more technical result.
This requires a bit of machinery to state, however when $L$ is given by a
smooth submanifold it reduces to the following.

\begin{definition}
  \label{LS smooth}Let $( M, \omega )$ be a symplectic 4--manifold with the
  decomposition $( L, E \rightarrow \Sigma )$ such that $L \hookrightarrow M$
  is a smooth Lagrangian submanifold of $M$.
  \begin{enumerate}
    \item Let $\mathcal{S}(M)$ denote the symplectomorphisms of $( M,
    \omega )$.

    \item Let $\mathcal{L}^{\tmop{sm}}$ denote the Lagrangian embeddings of
    $\phi\co L \hookrightarrow M$ which extend to symplectomorphisms of $M$.

    \item Let  $\mathcal{L}^{} \mathcal{E^{\tmop{sm}}}$ denote the space of
    pairs $( \psi, S )$ where $\psi \in \mathcal{L}^{\tmop{sm}}$ and $S$ is a
    symplectic embedded unparametrized surface which is abstractly
    symplectomorphic to $\Sigma$ and disjoint from $\psi ( L )$.

    \item Let $\mathcal{E}mb_{\omega} ( \Sigma, E )$ denote the
    space of unparametrized, embedded symplectic surfaces $S$ in $E
    \subset M$ such that $\omega [ S ] = \omega [ \Sigma ]$.
  \end{enumerate}
\end{definition}

\begin{theorem}[The main theorem when $L$ is smooth]
  \label{thm:spine smooth}Let $( M, \omega )$
  be a symplectic 4--manifold with Biran decomposition $( L, E \rightarrow
  \Sigma )$ such that $\phi\co L \hookrightarrow M$ is a smooth, Lagrangian
  submanifold. Then $\mathcal{S} ( M )$ is homotopy equivalent to
  $\mathcal{L}^{} \mathcal{E^{\tmop{sm}}}$.

  Moreover there is a fibration $\mathcal{L} \mathcal{E}^{\tmop{sm}}
  \rightarrow \mathcal{L^{\tmop{sm}}}$ whose fiber is homotopy equivalent to
  $\mathcal{E} m b_{\omega} ( \Sigma, E )$. When
  $\Sigma$ has genus 0, $\mathcal{E} m b_{\omega} ( \Sigma, E )$ is
  contractible and thus $\mathcal{S} ( M )$ is homotopy equivalent to
  $\mathcal{L}^{\tmop{sm}}$.
\end{theorem}

We note that in the case that $M$ is a surface the corresponding
symplectic embeddings $\mathcal{E} m b_{\omega} (
\Sigma, E )$ are just the embeddings of each point inside the
appropriate disc. Thus in the two-dimensional case $\mathcal{E} m
b_{\omega} ( \Sigma, E )$ is always contractible, and one recovers
(modulo concerns of $L$'s smoothness) Theorem \ref{thm:classical
sym}.  We show that in dimension 4, $\mathcal{E} m b_{\omega} (
\Sigma, E )$ is contractible if $\Sigma$ is a sphere. However the
proof relies heavily on the special properties of $J$--holomorphic
spheres in rational surfaces, and it is unclear how one may
generalize it to cases when $\Sigma$ has higher genus.

In general this theorem, and more generally Theorem
\ref{thm:Symp(M)is-homotopy-equivalent LSigma}, allow us to
separate the problem of understanding the topology of the
symplectomorphism group into two parts: embeddings of isotropic
skeletons up to symplectic equivalence, and $\mathcal{E} m
b_{\omega} ( \Sigma, E )$ the embeddings of symplectic surfaces in
a fixed homology class disjoint from the spine. This second
problem is universal, depending only on the genus, and self-intersection
of $\Sigma$ but not the ambient manifold $M$.

Understanding the higher homotopy groups of spaces of Lagrangian embeddings
$L \rightarrow M$ is quite hard.  Eliashberg and Polterovich showed that such
spaces are locally contractible, ie the space of Lagrangian embeddings of
$\mathbb{R}^2$ into $\mathbb{R}^4$ which are standard at $\infty$ is
contractible {\cite{PE2}}.  However, prior to this paper the author knows of
no computation when both domain and range are closed. Thus, while Theorem
\ref{thm:Symp(M)is-homotopy-equivalent LSigma} provides a satisfying
generalization it is, at present, difficult to use it to compute much about
the symplectomorphism group of a 4 manifold. However we can use it leverage
our knowledge of symplectomorphism groups into an understanding of spaces of
Lagrangian embeddings -- a result which should satisfy in proportion to our
previous frustration. We obtain the following corollaries in subsection
\ref{sec:embeddings}, showing that spaces of Lagrangian embeddings due indeed
have non-trivial global topology. Each is obtained by combining our result with
Gromov's computations of the symplectomorphism groups of $\mathbb{CP}^2$ and
$S^2 \times S^2$.

\begin{theorem}
  \label{lagrp2}The space of Lagrangian embeddings of $\mathbb{R
  P}^2 \hookrightarrow \mathbb{CP}^2$ isotopic to the standard one is
  homotopy equivalent to $\mathbb{PU} ( 3 )$
\end{theorem}

\begin{theorem}
  \label{laganti}Let $\omega$ be a symplectic form on $S^2 \times S^2$ such
  that $\omega [ S^2 \times pt ] = \omega [ pt \times S^2 ]$, then the space
  of Lagrangian embeddings $S^2 \hookrightarrow S^2 \times S^2$ isotopic to
  the anti-diagonal is homotopy equivalent to $SO ( 3 ) \times SO ( 3 )$.
\end{theorem}

Richard Hind has recently proven that every Lagrangian sphere in $S^2 \times
S^2$ (with the above symplectic structure) is isotopic to the anti-diagonal
{\cite{Hind}}. One suspects that his methods probably can be used to show that
every Lagrangian $\mathbb{RP}^2$ is isotopic to the standard one. Thus both
theorems above should actually compute the homotopy type of the full space of
Lagrangian embeddings.

Finally we note that the spine $L$ is smooth only in a few special, though
important, decompositions, and those which we apply to obtain the above two
corollaries may nearly exhaust them.  In general one requires the extra
generality and complexity of Theorem \ref{thm:Symp(M)is-homotopy-equivalent
LSigma}.  In particular, one should not hope to compute the space of
Lagrangian embeddings $L \rightarrow M$ for general domain and target via the
methods of this paper, without considerable and yet unseen adaptation.

\subsection{Acknowledgements}

This paper is a portion of the author's PhD thesis, produced under the
helpful guidance of his advisors Dusa McDuff and Dennis Sullivan. The author
would particularly like to thank Dusa McDuff for her tireless aid in preparing
this article for publication.  He would also like to thank the referee for his
helpful comments. {The author is supported in part by NSF grant DMS-0202553.}

\section{Definitions and precise statements}

\subsection{Kan complexes }

\subsubsection{Motivation: germs and Kan complexes} \label{sec:Kan}

In what follows we need to understand the space of symplectomorphisms fixing
an isotropic skeleton, \textit{and the germ of a neighborhood surrounding it.}
This forces us to work with ``spaces'' of germs of mappings. A germ does not
have a specified domain, and as a result most natural topologies on spaces of
germs have unwanted pathologies.

While the space of germs is difficult to define, compact families
of germs are transparent -- there is no difficulty with domains. In
this paper we are concerned here with weak homotopy type: the
study of compact families. Thus there is no real foundational
difficulty. We do however require some linguistic finesse, and Kan
complexes, simplicial sets that satisfy the extension condition,
provide exactly this.

One can find a straightforward introduction to Kan complexes in J.P. May's
book {\cite{May}}. However on a first reading the reader can safely ignore
these subtleties and replace each Kan complex with the appropriate ``space''.
Indeed, when the spine is a smooth Lagrangian manifold every symplectomorphism
fixing the spine also fixes its normal bundle, and thus the symplectomorphisms
fixing a neighborhood of the spine are a deformation retract of those fixing
the spine.

\subsubsection{A brief introduction to Kan complexes}

Since Kan complexes are perhaps unfamiliar to much of this paper's audience we
shall give a brief informal introduction here, based on that in May's book.

\paragraph{The basic example: singular simplices of a topological space}

Let $\Delta_n$ denote the $n$--simplex. If $X$ is a topological
space, the singular $n$--simplices $\Delta^n ( X )$ are the
continuous maps
\[ g\co \Delta_n \rightarrow X. \]
The singular simplices of $X$ are then given by $\Delta ( X ) =
\bigcup_{n \in \mathbb{N}} \Delta^n ( X )$.

Kan complexes are designed to axiomatize those properties of the simplices
$\Delta ( X )$ that are required to do homotopy theory -- to construct homotopy
groups, fibrations etc.  They will be particularly useful in our case
because as mentioned above there will be several ill-defined ``spaces'' in our
discussion whose {\tmem{simplices}} {\tmem{are well defined}}.

\paragraph{Degeneracy and face operators}

We begin with a graded set $K = \bigcup K_q$. The grading gives the
``dimension'' of the simplices. The first piece of extra structure we require
are the face and degeneracy operators $\partial_i\co K_q \rightarrow K_{q-1}$ and
$s_i\co K_q \rightarrow K_{q+1}$.  $\partial_i$ axiomatizes the notion of
passing from a singular simplex to its $i$th face. $s_i$ encodes
that of passing from an $n$--simplex to a degenerate $(n{+}1)$--simplex.  One
then has a list of axioms encoding the geometry of these operations
\[ \begin{array}{llll}
     \partial_i \partial_j & = & \partial_{j - 1} \partial_i & \tmop{if} i <
     j\\
     s_i s_j & = & s_{j + 1} s_i & \tmop{if} i \leq j\\
     \partial_i s_j & = & s_{j - 1} \partial_i & \tmop{if} i < j\\
     \partial_j s_j & = & \tmop{identity} & = \partial_{j + 1} s_j,\\
     \partial_i s_j & = & s_j \partial_{i - 1} & \tmop{if} i > j + 1
   \end{array} \]
A graded set with face and degeneracy operations is called a
{\tmdfn{simplicial set}}.

\paragraph{The extension condition}

The second necessary structure is that of {\tmem{extension}}. Suppose one has
a continuous map $f$ into $X$ which is defined on the union of $n$ faces (all
but one) of $\Delta_n$. Then since $\Delta_n$ retracts onto this union, one
can extend $f$ to be a singular $n$--simplex, a continuous map of all of
$\Delta_n$ into $X$.

One axiomatizes this by saying that if one has a collection of
$(n{-}1)$--simplices in $K_{n-1}$ whose degeneracy and face maps fit
together like this union of $n$ faces, then one can find an $n$--simplex
in $K_n$ filling them in.  The extension condition is also called the
Kan condition. It is required for homotopy of simplices to form an
equivalence relation.

\paragraph{Homotopy groups}

Given two singular simplices $f, g\co \Delta_n \rightarrow X$, sharing a
common boundary, we say that they are homotopic (relative to their boundary)
if there is an  $(n{+}1)$--simplex $h$ whose non-degenerate faces are $f$
and $g$.

\begin{figure}[ht!]
\begin{center}
\includegraphics{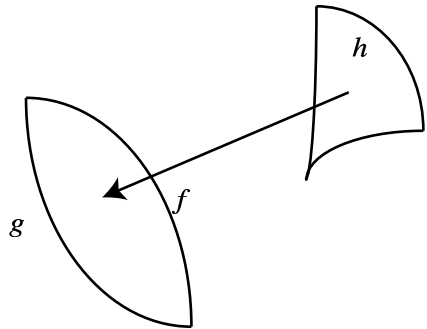}
\caption{The singular 1--simplices $f$ and $g$ are homotopic
because there is a 2--simplex $h$ making each a face.}
\label{figure:example}
\end{center}
\end{figure}

There is a natural generalization of this definition to  the homotopy of
simplices in Kan complexes.  Again one encodes the phenomena in terms of the
face and degeneracy operators. Homotopy of simplices is then an equivalence
relation.

Given a based topological space $( X, x_0 )$, one can think of $\pi_n ( X,
x_0 )$ as the homotopy classes of singular simplices $f\co \Delta_n \rightarrow
X$ such that the faces of $\Delta_n$ all map to $x_0$.  In this way one can
also define homotopy groups of a Kan complex.  One first defines a basepoint
by choosing an $x_0 \in K_0$, and considering it together with the sub complex
generated by its degeneracy maps.  This yields a subcomplex with one
completely degenerate simplex in each dimension that we denote also by $x_0$.
Then one can define $\pi_n ( K, x_0 )$ as the set of homotopy classes of
$n$--simplices in $K_n$ whose faces all lie in $x_0$. These have a group
structure when $n \geq 1$, which is Abelian for $n \geq 2$. One can also
define relative Kan homotopy groups, having the expected properties, in an
analogous manner.

\paragraph{Kan fibrations}

A Serre fibration of topological spaces is one that allows ``lifting'' of
homotopies of simplices.  A Kan fibration is defined analogously.  Just as for
spaces, a Kan fibration of based complexes induces a corresponding long exact
sequence of homotopy groups,

\paragraph{The geometric realization}

Given a Kan complex $K$ one can, in a natural way, construct a topological
space $|K|$ by constructing  a topological simplex for each simplex in $K$ and
then ensuring that the faces and and degenerations act appropriately:
\[ |K| = \bigcup_{n \in \mathbb{N}} K_n \times \Delta_n / \sim \]
Where $\sim$ is an equivalence relation given by
\begin{eqnarray*}
  ( \partial_i k_n, u_{n - 1} ) & \sim & ( k_n, \partial_i u_{n - 1} ),\\
  ( s_i k_n, u_{n + 1} ) & \sim & ( k_n, s_i u_{n + 1} ),
\end{eqnarray*}
for  $k_n \in K_n, u_{n-1} \in \Delta_{n-1}$, and $u_{n+1} \in
\Delta_{n+1}$.
$\Delta(|K|)$ is  weakly (Kan) homotopy equivalent to $K$. Similarly,
for any topological space $X$, $| \Delta ( X ) |$ is weakly homotopy
equivalent to $X$.

\subsection{Biran decompositions}

Let $( M, \omega )$ denote a symplectic $4$--manifold. For any symplectic
manifold $( N, \sigma )$, \textbf{$\mathcal{S} ( N, \sigma )$} denotes the
symplectomorphisms of $( N, \sigma )$. If either $N$ or $\sigma$ is clear from
the context they will be omitted. \textbf{$\mathcal{D} ( M )$} denotes the
diffeomorphisms of $M$. Again, if $M$ is clear we will omit it.

\begin{definition}
  A smoothly embedded cell complex consists of:
  \begin{enumerate}
    \item An abstract smooth cell complex $C$ in which the interior of each
    cell is endowed with a smooth structure.

    \item A continuous map
    \[ i\co C \hookrightarrow M \]
    which is a smooth embedding when restricted to the interior of each cell
    in $C.$
  \end{enumerate}
  We say that a smoothly embedded cell complex is {\tmdfn{isotropic}} with
  respect to a symplectic structure $\omega$, if $i^{\ast} ( \omega ) = 0$ on
  the interior of each cell.
\end{definition}

\begin{definition}
  Let $( M, \omega )$ be a symplectic manifold. Let $J$ be an almost complex
  structure compatible with $\omega .$ Let $\Sigma_{\lambda}$ a symplectic
  hypersurface of $M$, Poincare dual to $\lambda [ \omega ]$, and such that:
  \begin{enumerate}
    \item There is a smoothly embedded, isotropic cell complex $L_{\lambda}$
    disjoint from $\Sigma_{\lambda}$. In what follows we will call this cell
    complex an \emph{isotropic skeleton} of $M$.

    \item $M - L_{\lambda}$ is an open symplectic disc bundle $E$ over
    $\Sigma_{\lambda}$, such that the fibers have area $\frac{1}{\lambda}$
    with respect to $\omega$. This bundle is symplectomorphic to the unit disc
    bundle in the normal bundle to $\Sigma_{\lambda}$ with symplectic form
    \[ \pi^{\ast} \omega |_{\Sigma} + d ( r^2 \alpha )/\lambda, \]
    where $r$ is the radial coordinate in the fiber, $\alpha$ is the
    connection 1--form coming from the hermitian metric $\omega ( \cdot, J
    \cdot )$ on the normal bundle, and $\alpha$ is normalized so that its total
    integral around the boundary of a fiber is $\frac{1}{\lambda}$.
  \end{enumerate}
  We call such a configuration $( L_{\lambda}, E \rightarrow \Sigma_{\lambda}
  )$ a \emph{decomposition} of $M$.
\end{definition}

\begin{theorem}[Biran \cite{biran:Lspine}]
  \label{thm:(Biran)-Let-M} If M is a
  smooth, projective variety, then there is a decomposition $( L_{\lambda}, E
  \rightarrow \Sigma_{\lambda} )$.
\end{theorem}

Biran conjectures that his proofs generalize to any symplectic manifold
whose symplectic structure $\omega$ has a rational cohomology class.
This would rely on Donaldson's result {\cite{Donaldson}} implying that
in such a setting one has a symplectic hypersurface Poincare dual to $[
\lambda \omega]$.

\subsection{Notation and a precise statement of our main theorem}

Let $( L, E \rightarrow \Sigma )$ be a decomposition of the
symplectic $4$--manifold $( M, \omega )$. Denote the symplectic
embeddings of $\Sigma$ into $M$ by $\mathcal{E} m b_{\omega} (
\Sigma, M )$.

\begin{definition}
  Denote by $\mathcal{L}$ the Kan complex of isotropic embeddings of $L$ into
  $M$.

  An  $n$--simplex in $\mathcal{L}$ is an equivalence class given
  by:
  \begin{enumerate}
    \item A neighborhood $U$ of $L$ in $M$.

    \item A continuous map $\phi\co \Delta_n \rightarrow \mathcal{S} ( U, M
    )$, where $\mathcal{S} ( U, M )$ denotes the symplectic embeddings of $U$
    into $M$, which admits an extension to a symplectomorphism of all of $M$.
  \end{enumerate}
  Two such pairs $( U_1, \phi_1 )$ and $( U_2, \phi_2 )$ are equivalent if
  there exists a neighborhood $U_3$ of $L$ such that $U_3 \subset U_1$, $U_3
  \subset U_2$ and
  \[ \phi_1 |_{U_3} = \phi_2 |_{U_3} \]
  The face maps $\partial_i$ are given by restricting $\phi$ to the faces of
  $\Delta_n$. The degeneration maps $s_i$ are given by pre-composing $\phi$
  with the degeneration maps of $\Delta_n$.
\end{definition}

\begin{definition}
  Denote by $\mathcal{L} \mathcal{E}$ \textbf{}the following Kan complex:

  an $n$--simplex consists of a pair $( \text{$( U, \phi )$}, \psi )$ where:
  \begin{enumerate}
    \item $( U, \phi )$ is an $n$--simplex in $\mathcal{L}$.

    \item A continuous map $\psi\co \Delta_n \rightarrow \mathcal{E} m b_{\omega} ( \Sigma, M )$ such that for each $x \in \Delta_n$, $\psi ( x
    )$ is disjoint from $\phi ( L )$.
  \end{enumerate}
  The face maps $\partial_i$ are given by restricting $\phi$ and $\psi$
  to the faces of
  $\Delta_n$. The degeneration maps $s_i$ are given by pre-composing
  $\phi$ and $\psi$ with the degeneration maps of $\Delta_n$.
\end{definition}

Note that $\mathcal{S} ( M )$ acts on $\mathcal{L} \mathcal{E}$. For any
symplectomorphism will carry any simplex of skeletons  to another , and will
preserves the the homology class of $\Sigma$ as it is Poincare dual to
$\lambda [ \omega ]$.

This paper is devoted to the proof and application of the following theorem.

\begin{theorem}
  \label{thm:Symp(M)is-homotopy-equivalent LSigma}Let $( M, \omega )$ be a
  symplectic $4 -$manifold with decomposition $( L, E \rightarrow \Sigma
  )$. Then the orbit map gives a weak homotopy equivalence of  $\mathcal{S}( M
  )$ with $\mathcal{|L} \mathcal{E} |$. Moreover there is a fibration
  $\mathcal{|L} \mathcal{E} | \rightarrow | \mathcal{L|}$ whose fiber is given
  by $\mathcal{E} m b_{\omega} ( \Sigma, E )$, the
  symplectic embeddings of $\Sigma$ into $E$. When $\Sigma$ has genus 0,
  $\mathcal{E} m b_{\omega} ( \Sigma, E )$ is contractible.
\end{theorem}

Here and elsewhere, when we say a map $f\co X \rightarrow Y$ is a homotopy
equivalence we mean that:
\begin{enumerate}
  \item $f$ gives a bijection between the connected components.

  \item If one chooses a base point $x_i$ in each connected component $X_i$ of
  $X$, and chooses $f ( x_i )$ to be the basepoint of corresponding connected
  component $Y_i$ of $Y$; then the resulting map $f\co ( X_i, x_i ) \rightarrow
  ( Y_i, f ( x_i ) )$ gives a weak homotopy equivalence.
\end{enumerate}
We do not require that $M$ be K\"ahler, only that it has a decomposition.
However we do restrict ourselves to dimension 4, as we rely heavily on the
properties of $J$--holomorphic spheres in this dimension.

\subsection{Conventions}

Note that both the symplectomorphisms $\mathcal{S} ( M, \omega )$ and the
space of pairs $\mathcal{L} \mathcal{E}$ are invariant under scaling the
symplectic structure by a constant factor. Thus we safely replace $\omega$ by
$\lambda \omega$ and reduce to the case that the class of $\Sigma$ and the
symplectic form are Poincare dual, and $E$'s fibers have area 1.

\paragraph{Groups and actions}

Throughout this paper we will be computing and comparing the stabilizers of
the action of various groups on various geometric objects. To keep our heads
straight it will be helpful to adopt a few notational guidelines.
\begin{enumerate}
  \item $G_S$ denotes the elements in the group $G$ which preserve the set
  $S$. That is, $G_S = \{ g \in G : g ( S ) \subset S \}$.

  \item $G_{\mathbb{S}}$ denotes the elements in the group $G$ which fix the
  set $S$. That is, $G_\mathbb{S} = \{ g \in G : g ( s ) = s, \forall s \in S \}$.  This notation
  is required only in the next subsection.

  \item $G_{\overline{S}}$ denotes the elements in the group
  $G$ which fix
  both a set S and a framing of that set. These are
  $$\{ g \in G : \exists \textrm{ a \tmop{neighborhood }} N_S \supset S :
  g ( s ) = s, \forall s \in N_S \}.$$
  We endow $G_{\overline{S^{}}}$ with the direct limit topology and
  say that these are the elements of $G$ which fix a framing of $S$.

  \item If we fix or preserve more than one set we will denote this by
  separating the two with a comma. For example: $G_{X, \widetilde{Y}}$
  denotes the elements in $G$ which preserve $X$ and fix both $Y$ and
  a framing.
\end{enumerate}

\subsection{\label{simpler smooth}Simplifications when $L$ is smooth: Proof
of Theorem \ref{thm:spine smooth}}

\begin{proof}
  If the skeleton of the decomposition $L$ is a smooth, Lagrangian,
  submanifold one does not have to introduce Kan complexes. In this case
  $\mathcal{S} ( M )$ is homotopy equivalent to $\mathcal{L}^{}
  \mathcal{E^{\tmop{sm}}}$ (see Definition \ref{LS smooth}).

  The difference between this result (Theorem \ref{thm:spine smooth}) and the
  more general statement of Theorem \ref{thm:Symp(M)is-homotopy-equivalent
  LSigma} is as follows:  Theorem \ref{thm:Symp(M)is-homotopy-equivalent
  LSigma} states that $\mathcal{S}_{\overline{L}, \Sigma}$, the symplectomorphisms
  which fix a framing of $L$ and preserve $\Sigma$, are contractible,
  while Theorem \ref{thm:spine smooth} requires that
  $\mathcal{S}_{\mathbb{L}, \Sigma}$, the symplectomorphisms which
  only fix $L$ and preserve $\Sigma$, are contractible.

  We will now show that the natural map  $\mathcal{S}_{\overline{L}, \Sigma}
  \rightarrow \mathcal{S}_{\mathbb{L}, \Sigma}$ is a homotopy equivalence,
  and thus that Theorem \ref{thm:spine smooth} follows from Theorem
  \ref{thm:Symp(M)is-homotopy-equivalent LSigma}. If $\phi \in
  \mathcal{L^{\tmop{sm}}}$, every extension of $\phi$ to a symplectomorphism
  of $M$ induces the same map on the normal bundle $N_L$ of $L$. This is an
  immediate consequence of the corresponding linear statement: If $( V,
  \omega )$ is a symplectic vector space  with  Lagrangian subspace L,
  $\phi_1$ and $\phi_2$ are linear symplectic automorphisms of $V$ such that
  $\phi_1 |_L = \phi_2 |_L$ then each induces the same map $V / L \rightarrow
  V / \phi_i ( L )$.

  Thus Weinstein's neighborhood theorem (in parameters) implies that the map
  $\mathcal{L} \hookrightarrow \Delta(\mathcal{L}^{\tmop{sm}})$ is a weak
  homotopy equivalence. The actions of
  $\mathcal{S} ( M )$ on each space yield a morphism of fibrations:
  \[ \begin{CD}
  \Delta(\mathcal{S}_{\overline{L}} ) @>i_1>>
    \Delta(\text{$\mathcal{S}_{\mathbb{L}}$}) \\
  @VVV    @VVV \\
  \Delta(\mathcal{S}) @>(\text{id})>> \Delta(\mathcal{S}) \\
  @VVV    @VVV \\
  \mathcal{L} @>i_2>> \Delta(\mathcal{L}^{\tmop{sm})}
  \end{CD} \]

  By applying the 5--Lemma (Lemma \ref{lem:5Lemma}) to the associated long
  exact sequence of Kan homotopy groups we see that
  $\Delta(\mathcal{S}_{\overline{L}}) \overset{i_1}{\hookrightarrow}
  \Delta(\mathcal{S}_\mathbb{L})$ must also be a homotopy equivalence, so
  $\mathcal{S}_{\overline{L}} \rightarrow S_{\mathbb{L}}$ and
  $\mathcal{S}_{\overline{L}, \Sigma} \rightarrow \mathcal{S}_{\mathbb{L},
  \Sigma}$ are also homotopy equivalences.
\end{proof}

\subsection{Applications to spaces of Lagrangian embeddings
\label{sec:embeddings}}

We now apply Theorem \ref{thm:spine smooth} to compute spaces of Lagrangian
submanifolds in cases where we know the homotopy type of $\mathcal{S} ( M )$.
We show that:
\begin{enumerate}
  \item  The{\tmem{}} space of Lagrangian embeddings of $\mathbb{R P}^2
  \hookrightarrow \mathbb{CP}^2$ isotopic to the standard one is homotopy
  equivalent to $\mathbb{PU} ( 3 )$.

  \item  If $\omega$ is a symplectic form on $S^2 \times S^2$ such that
  $\omega [ S^2 \times pt ] = \omega [ pt \times S^2 ]$ then the space of
  Lagrangian embeddings $S^2 \hookrightarrow S^2 \times S^2$ isotopic to the
  anti-diagonal is homotopy equivalent to $\tmop{SO} ( 3 ) \times \tmop{SO} (
  3 )$.
\end{enumerate}
These are Theorems \ref{lagrp2} and \ref{laganti} respectively.  In
{\cite{biran:Lspine}} Biran computes decompositions of  $\mathbb{CP}^2$ and
$( S^2 \times S^2, \omega )$ with a symplectic structure such that $\omega [
S^2 \times pt ] = \omega [ pt \times S^2 ]$.

\begin{proposition}[Biran]$\phantom{99}$
  \begin{enumerate}
    \item $\mathbb{CP}^2$ admits a decomposition where $\Sigma$ is a
    quadric (and
    thus a sphere) and $L$ the standard $\mathbb{RP}^2 \hookrightarrow
    \mathbb{CP}^2$.

    \item $( S^2 \times S^2, \omega )$ with a symplectic structure such that
    $\omega [ S^2 \times pt ] = \omega [ pt \times S^2 ]$ admits a
    decomposition with $\Sigma$ the diagonal and $L$ the anti-diagonal.
  \end{enumerate}
\end{proposition}

In each case the decomposition has $L$ a smooth manifold with $H_1
( L )$ torsion, and $\Sigma$ a sphere. Denote the identity
component of $\mathcal{S} ( M )$ by $\mathcal{S} ( M )^{\circ}$.
Similarly we denote the identity component of
$\mathcal{L}^{\tmop{sm}}$ by $\mathcal{L}_{\circ}^{\tmop{sm}}$. In
the two decompositions above $\mathcal{L_{\circ}^{\tmop{sm}}}$ is
the space of all Lagrangian embeddings of $L$ isotopic to the
original. For in each case $H_1 ( L )$ is torsion $H_2 ( M )
\otimes \mathbb{R} \rightarrow H_2 ( M, L ) \otimes \mathbb{R}$ is
surjective, and thus every isotopy of $\phi$ can be induced by an
ambient symplectic isotopy of $M$.

The orbit map
\[ \text{$\mathcal{S} ( M )^{\circ}$} \rightarrow
   \mathcal{L}_{\circ}^{\tmop{sm}} \]
is a surjective fibration. By Theorem \ref{thm:spine smooth} its
fiber is homotopy equivalent to  $\mathcal{E} m b_{\omega} (
\Sigma, E )$. Moreover since $\Sigma$ is a sphere, Theorem
\ref{thm:spine smooth} states that $\mathcal{E} m b_{\omega} (
\Sigma, E )$ is contractible. Thus the map is a homotopy
equivalence.

Finally, since the identity components of $\mathcal{S} ( \mathbb{CP}^2 )$
and $\mathcal{S} ( ( S^2 \times S^2, \omega )$, are homotopy equivalent to
$\mathbb{P U} ( 3 )$ and $\tmop{SO} ( 3 ) \times \tmop{SO} ( 3 )$
respectively {\cite{phol}} Theorems \ref{lagrp2} and \ref{laganti} follow. In
general we have the following corollary:

\begin{corollary}
  \label{cor:sphere makes symp same as lagrange}Let $( M, \omega )$ be a
  symplectic 4--manifold with the decomposition $( L, E \rightarrow \Sigma )$
  such that $\phi\co L \hookrightarrow M$ is a smooth submanifold of $M$, $H_2
  ( M ) \otimes \mathbb{R} \rightarrow H_2 ( M, L ) \otimes \mathbb{R}$ is
  surjective, and $\Sigma$ is a sphere. Then the space $\mathcal{L}_{\phi}$ of
  Lagrangian embeddings isotopic to $\phi$ is homotopy equivalent to the
  identity component of $\mathcal{S} ( M )$.
\end{corollary}

\section{Reduction to a problem on rational surfaces} \label{sec:reduction to
rational}

We seek in this section to perform a fiber-wise compactification of (most of)
the disc bundle $E$ into a symplectic sphere bundle  $\widehat{E}_{\epsilon}
\rightarrow \Sigma$. This sphere bundle will have two distinguished symplectic
sections $Z_0$, which is identified with $\Sigma$ under the compactification,
and $Z_{\infty}$, which is the image of a neighborhood of $L$.

Then our theorem is reduced to a problem on rational surfaces:

\begin{definition}
  \label{def:cohomology assumption} Let $( X, \omega )$ be a symplectic
  4--manifold which is a symplectic sphere bundle over a surface $\Sigma$ with
  a symplectic form $\omega$. Let $Z_0$ be a symplectic section of $F$ such
  that $[ Z_0 ] \cdot [ Z_0 ] = k \geq 0$.  Denote the homology class of the
  fiber by $[ F ]$. We say that $( X, Z_0, F, \omega )$  {\emph{satisfies
  the cohomology assumption}} if $[ \omega ] = a \cdot PD ( [ Z_0 ] ) + b
  \cdot PD ( [ F ] )$ with $a, b > 0$ both positive real numbers.

\end{definition}

\begin{proposition}[Problem on rational surfaces]
\label{rational surface}
  Let $( X, Z_0, F,
  \omega )$ be a symplectic 4--manifold  admitting a symplectic fibration by
  2--spheres by $F$,  $Z_0$ be a symplectic section of $F$ such that $[ Z_0 ]
  \cdot [ Z_0 ] = k \geq 0$. Suppose that $( X, Z_0, F, \omega )$ satisfies
  the cohomology assumption.  Let $Z_{\infty}$ be a symplectic section of $F$
  which is disjoint from $F$, and thus has self-intersection $[ Z_{\infty} ]
  \cdot [ Z_{\infty} ] = - k$.

  Then the space of symplectomorphisms $\mathcal{S^{} ( \text{}} X
  )_{\overline{\infty}}$ of $X$ which fix a framing of $Z_{\infty}$ acts
  transitively on the space $\mathcal{E}_{\infty}$ of unparametrized,embedded
  symplectic surfaces $S$ of $X$ disjoint from $Z_{\infty}$ with contractible
  stabilizer. Moreover, $\mathcal{E}_{\infty}$ is contractible when $\Sigma$
  has genus 0.
\end{proposition}

The compactification $\widehat{E}_{\epsilon}$ serves two roles: it
allows us to play in the more comfortable compact terrain, and
it converts the problem of computing the stabilizer of an
isotropic object (the spine $L$) to that of a symplectic object
(the section $Z_{\infty}$). For this dual service we pay a price:
we cannot compactify all of $M - L$, and must be satisfied with
compactifying the complement of neighborhood of $L$.

The compactification is described in \ref{sub:Compactify}, the reduction of
our Theorem to Proposition \ref{rational surface}  is described in
\ref{sub:Restatement-of-Theorem}. Its proof is described in section
\ref{sec:ContractibilityandTransitivity} once we have developed the necessary
background in sections \ref{sec:Curves-and-Fibrations} and 5.

\subsection{Compactification of $E_{1 - \epsilon}$ via symplectic cutting (a
la Lerman)} \label{sub:Compactify}

We apply the techniques of {\cite{Lerman}} to the present situation. Consider
$E \subset M$ as the open unit disc bundle in the normal bundle $N_{\Sigma}$
to $\Sigma$. Denote by $E_{1 - \epsilon} \subset E$ the set $\{ x \in E :
||x|| \leq 1 - \epsilon \}$.

\begin{lemma}
  \label{lem:Compactification}There is a surjective $C^{\infty}$ map $\Psi\co
  E_{1 - \epsilon} \rightarrow \widehat{E}_{\epsilon}$ where
  $\widehat{E}_{\epsilon}$ is a symplectic sphere bundle over $\Sigma$.
  Topologically, $\Psi$ is given by the collapse of the boundary circle
  in each fiber of the disc bundle $E_{1 - \varepsilon}$.
  \begin{enumerate}
    \item $\Psi$ is a symplectomorphism on the interior of $E_{1 - \epsilon}$.

    \item $\Psi$ maps the boundary of $E_{1 - \epsilon}$ to a symplectic
    section $Z_{\infty}$ of this bundle whose self-intersection is $- k$.

    \item $\Psi$ maps the zero section {\tmem{{\tmem{}}}}of $E_{1 - \epsilon}$
    to a symplectic section $Z_0$ whose self-intersection is $k$.
    Moreover the symplectic form $\omega$ on $\widehat{E}_{\epsilon}$ has
    cohomology class $( 1 - \epsilon ) PD ( [ Z_0 ] ) + \epsilon kPD ( [ F ]
    )$ where $[ F ]$denotes the class of the fiber of
    $\widehat{E}_{\varepsilon}$. Thus $\widehat{E}_{\varepsilon}$
    satisfies the cohomology assumption (Definition \ref{def:cohomology
    assumption}).
  \end{enumerate}
\end{lemma}

\begin{proof}
  The bundle $E \rightarrow \Sigma$ is given as the unit disc bundle in
  $N_{\Sigma}$ in the hermitian metric induced from that on $M$. Place the
  coordinates $(r,t)$ on the fiber of $E = D^2$, where $r$ is a radial
  coordinate $r = |w|$, and the angular coordinate $t$ lies in $[ 0, 1 ]$,
  (ie $t = \frac{\theta}{2 \pi}$) Then the symplectic structure on E is
  given by
  \[ \pi^{\ast} \omega |_{\Sigma} + d ( |w|^2 \alpha ), \]
  where $\alpha$ is a connection one form. This structure is invariant under
  the circle action $S(t)$ given by the Hamiltonian function $\mu = |w|^2$.

  We now consider the $S^1$ action $P ( t )$ on the product
  \[ ( E \times C, \omega \oplus \tau ), \]
  where $C$ denotes the complex numbers and $\tau$ denotes their standard
  complex structure, scaled by the constant factor $\frac{1}{\pi}$. The action
  is given by
  \[ P ( t ) ( m, z ) = ( S ( t ) m, e^{2 \pi it} z ), \]
  where $P ( t )$ is Hamiltonian with function
  \[ \zeta = \mu + ||z||^2. \]
  Let $\widehat{E}$ be the symplectic reduction of $( ( E \times C, \omega
  \oplus \tau ), P ( t ) )$ along the level set $\zeta = 1 - \epsilon$. The
  level set
  \[ \zeta_{1 - \epsilon} := \{ ( m, z ) : \zeta ( m, z ) = 1 - \epsilon \}
  \]
  is the disjoint union of
  \[ \left\{ ( m, z ) : \mu ( m ) < 1 - \epsilon
     \textrm{,} z = e^{2 \pi it} \sqrt{\mu ( m ) - ( 1 - \epsilon )} \right\}
     \] 
  $$\left\{ ( m, 0 ) : \mu ( m ) = 1 - \epsilon \right\}\leqno{\rm and}$$
  where both members of the disjoint union are invariant under the $S^1$
  action. The map $i\co E_{1 - \epsilon} \rightarrow E \times C$ given by
  \[ i( m ) = ( m, \sqrt{( 1 - \epsilon ) - \mu ( m )} \]
  is a symplectic embedding, whose image is contained in the level set
  $\zeta_{1 - \epsilon}$. I claim that the composition of $i$ with the
  quotient
  \[ \pi_Q\co\zeta_{1 - \epsilon} \rightarrow \zeta_{1 - \epsilon} / S^1 =
     \widehat{E}_{\epsilon} \]
  of $\zeta_{1 - \epsilon}$ by $P ( t )$
  gives a map
  \[ \psi = \pi_Q \circ i\co E_{1 - \epsilon} \rightarrow
     \widehat{E}_{\epsilon} \]
  with the properties above.
\end{proof}

\begin{remark}
  Symplectic structures on rational surfaces are classified up to
  diffeomorphism by their cohomology class {\cite{Lalonde}}. Thus
  condition (4) determines the symplectic structure.
\end{remark}

\subsection{\label{sub:Restatement-of-Theorem}Translation of the main theorem to
compactification}

This subsection is devoted to carrying out the reduction to rational surfaces
described in the introduction to section \ref{sec:reduction to rational}.

\paragraph{$\mathcal{S}_{\overline{L^{}}}$ acts transitively
on $\mathcal{E} m b_{\omega} ( \Sigma, E )$ with
contractible stabilizer}\nl I claim that the natural map
\[ \phi_{( \Sigma )}\co \mathcal{S}_{\widetilde{L^{}}} \rightarrow
   \text{$\mathcal{E} m b_{\omega} ( \Sigma, E )$} \]
is a surjective homotopy equivalence.

Denote by $\mathcal{E}_{\varepsilon}$ the space of
unparametrized, embedded symplectic surfaces $S$ in
$E_{1 - \epsilon} \subset M$ such that $\omega [ S ] = \omega [
\Sigma ]$, and denote by$ \mathcal{S}_{\overline{\varepsilon}}$ the
symplectomorphisms of $M$ which fix $M \backslash
E_{1 - \varepsilon}$ (a neighborhood of the isotropic skeleton
$L$).  Now consider the action of each
$\mathcal{S}_{\overline{\varepsilon}}$ on $\mathcal{E}_{\varepsilon}$.
The resulting orbit maps $\phi_i$ yield a morphism of direct
systems
\[ \begin{array}{ccccccc}
     \mathcal{S}_{\overline{\varepsilon_1}} & \hookrightarrow &
     \mathcal{S}_{\overline{\varepsilon_2}} & \hookrightarrow & \cdots &
     \hookrightarrow & \mathcal{S}_{\overline{L^{}}}\\
     \Big\downarrow \rlap{$\phi_1$} &  & \Big\downarrow\rlap{$\phi_2$} &
      & \cdots &  & \Big\downarrow \rlap{$\phi$}\\
     \mathcal{E}_{\varepsilon_1} & \hookrightarrow &
     \mathcal{E}_{\varepsilon_2} & \hookrightarrow & \cdots & \hookrightarrow &
     \text{ $\mathcal{E} m b_{\omega} ( \Sigma, E )$}
  \end{array} \]
for $\epsilon_i > \epsilon_{i + 1} > 0$.
I claim that each $\phi_i$ is a surjective homotopy equivalence, and thus
$\phi$ is as well. Note that the compactification $\Psi$  induces:
\begin{enumerate}
  \item A homeomorphism $\Psi_{\Sigma}\co \mathcal{E}_{\varepsilon} \rightarrow
  \mathcal{E}^{\epsilon}_{\infty}$, the space of unparametrized, embedded
  symplectic surfaces $S$ in $\widehat{E}_{\epsilon}$ disjoint from
  $Z_{\infty}$.

  \item A homeomorphism $\Psi_S\co \mathcal{S}_{\overline{\varepsilon}} \to
  \mathcal{S^{} ( \text{}} \widehat{E}_{\epsilon} )_{\overline{\infty}}$, the
  space of symplectomorphisms  $\widehat{E}_{\epsilon}$ which fix a framing
  of $Z_{\infty}$.
\end{enumerate}
Suppose we have proven Proposition \ref{rational
surface}. Then $\widehat{E}_{\epsilon}$ satisfies the cohomology
assumption and thus the map  $\mathcal{S^{} ( \text{}}
\widehat{E}_{\epsilon} )_{\overline{\infty}} \rightarrow
\mathcal{E}^{\varepsilon}_{\infty}$ is a surjective fibration with
contractible fiber. Thus we have
\[ \mathcal{S}_{\overline{\varepsilon^{}}} \leftrightarrow \mathcal{S^{} (
   \text{}} \widehat{E}_{\epsilon} )_{\overline{\infty}}
   \longrightarrow%_{\textrm{( \tmop{Proposition} \ref{rational surface} )}}
   \mathcal{E}^{\epsilon}_{\infty} \leftrightarrow \mathcal{E}_{\varepsilon},
\]
where the second map is given by Proposition~\ref{rational surface},
and  $\phi_{( \Sigma )}\co \mathcal{S}_{\overline{\varepsilon}} \rightarrow
\mathcal{E}_{\varepsilon}$ is a (weak) surjective homotopy equivalence for
$0 < \epsilon < 1$.  So  $\phi$ must also be a (weak) surjective
homotopy equivalence.

\paragraph{$\mathcal{S} ( M ) \simeq \mathcal{|L^{} E|}$} We consider the
action of $\mathcal{S} ( M )$ on $\mathcal{L} \mathcal{E}$. I claim that
the map $\phi_{( L, \Sigma )}\co \Delta ( \mathcal{S} ) \rightarrow
\mathcal{L^{} E}$
is a homotopy equivalence.  This will imply that the associated map on
geometric realizations $|\phi_{( L, \Sigma )}|\co |\Delta(\mathcal{S})|
\rightarrow \mathcal{|L^{} E|}$ is also a homotopy equivalence.
For since the natural inclusion $S \hookrightarrow | \Delta ( S ) |$
is always a homotopy equivalence, so is the composition
\[ S \hookrightarrow | \Delta ( \mathcal{S} ) | \rightarrow |\mathcal{LE}|. \]
This will prove our main theorem (Theorem~\ref{thm:Symp(M)is-homotopy-equivalent
LSigma}).

Consider then the Kan fibration
\[ \pi\co \Delta ( \mathcal{S} ) \rightarrow \mathcal{L E} \]
where $\mathcal{S}_{\overline{L^{}}, \Sigma}$ denotes the
symplectomorphisms of $M$ which fix a framing of $L$ and preserve
$\Sigma$. $\mathcal{L}$ is by definition a homogeneous space for
$\mathcal{S}$. That $\mathcal{L E}$ is as well follows from the
transitive action of $\mathcal{S}_{\overline{L^{}}}$ on
$\mathcal{E} m b_{\omega} ( \Sigma, E )$ proved above. Thus it is
enough to prove that $\phi_{( L, \Sigma )}$ is a homotopy
equivalence when the fibration is restricted to any connected
component.   The fiber of $\pi$ is then homotopy equivalent to
$\mathcal{S}_{\overline{L^{}}, \Sigma}$, and thus contractible.

\paragraph{If $\Sigma$ is a sphere, $ \mathcal{E} m b_{\omega} ( \Sigma, E )$
is contractible} $\mathcal{E} m b_{\omega} ( \Sigma, E )$ is given
by the direct limit
\[ \mathcal{E}_{\varepsilon_1} \hookrightarrow \mathcal{E}_{\varepsilon_2}
   \hookrightarrow \cdots \hookrightarrow \mathcal{E}mb_{\omega}(\Sigma,E) \]
as each $\mathcal{E}_{\varepsilon i} \simeq 
\mathcal{E}^{\epsilon_i}_{\infty}$ is thus contractible, so must
$\mathcal{E} m b_{\omega} ( \Sigma, E )$ be contractible.

\section{A softer symplectomorphism group of a rational surface}
\label{sec:Curves-and-Fibrations}
\begin{definition}
  Let  $( X, \omega )$ be a symplectic 4--manifold  admitting a
  symplectic fibration by 2--spheres by $F$, which satisfies the cohomology
  assumption. Let $Z_0$ be a symplectic section of $F$ such that $[ Z_0 ]
  \cdot [ Z_0 ] = k \geq 0$.  Let $Z_{\infty}$ be a symplectic section of $F$
  which is disjoint from $F$, and thus has self-intersection $[ Z_{\infty} ]
  \cdot [ Z_{\infty} ] = - k$.
\end{definition}

We will now construct a large open neighborhood of the symplectomorphisms
$\mathcal{S (} X )$ within the diffeomorphism group $\mathcal{D} ( X )$. This
neighborhood will have the same homotopy type as $\mathcal{S (} X )$, but it
will be far easier to work with. In particular, it will be much easier to
understand the ``action'' of this neighborhood on various objects. Various
such ``softenings'' have been used throughout the study the homotopy type of
symplectomorphism groups,  beginning with Gromov's initial work
{\cite{phol}} and continuing with the work of Abreu and McDuff
{\cite{abreu1,abreu}}.  The particular neighborhood we construct has its roots
in McDuff's inflation argument , used previously to classify the symplectic
structures on rational surfaces {\cite{Dusa:inflation}}.

We abbreviate $\mathcal{S} ( X )$ by $\mathcal{S} $ in this section and the
next.

\begin{definition}
  Denote by $\mathcal{\mathcal{F Z}}^{}$ the space of all triples $( F_S,
  S_0, S_{\infty} )$ where $F_S$ is a symplectic fibration of $X$ by two
  spheres in class $[ F ]$  and $S_0$ and $S_{\infty}$ are symplectic
  sections such that $[ S_0 ] = [ Z_0 ]$ and $[ S_{\infty} ] = [ Z_{\infty}
  ]$.
\end{definition}

\begin{definition}
  Denote by $\mathcal{S}^{\tmop{soft}}$ the diffeomorphisms $g$ of $X$ which
  preserve $H^2$ and such that the triple $( g ( F ), g ( Z_0 ), g (
  Z_{\infty} ) ) \in \text{$\mathcal{\mathcal{F Z}}^{}$} .$
\end{definition}

\begin{proposition}
  \label{S in Soft a retract} Let $Z_0$ be a symplectic section of $F$ such
  that $[ Z_0 ] \cdot [ Z_0 ] = k \geq 0$, and let $Z_{\infty}$ be a
  symplectic section of $F$ which is disjoint from $Z_0$. Then the inclusion
  $i\co\mathcal{S (} X ) \hookrightarrow \mathcal{S}^{soft} $ is a
  weak deformation retract.
\end{proposition} \label{Symp a retract of comp Diff}

\proof \label{lem:symp acts trivially on H}
  First we show that $\mathcal{S} \subset \mathcal{S}^{\tmop{soft}}$.
  If $\gamma \in \mathcal{S}$ it is
  clear that each member of the triple
  $$( \gamma ( F ), \gamma ( Z_0, ) \gamma ( Z_{\infty} ) )$$
  is symplectic. What is required then is to show
  that each such $\gamma$ preserves $H_2$. $\gamma$ preserves $\omega$, and
  thus also the cohomology class
  \[ [ \omega ] = ( 1 - \epsilon ) PD ( [ Z_0 ] ) + \epsilon kPD ( [ F ]). \]
  Thus, as $[ \omega ]$ is Poincare dual to
  \[ ( 1 - \epsilon ) ( [ Z_0 ] ) + \epsilon k ( [ F ] ), \]
  $\gamma$ preserves this homology class as well. As $[ Z_0 ]$ and $[ F_Z ]$
  together span $H_2 ( \widehat{E}_{\epsilon} )$, it is enough for us to
  show that $\gamma$ preserves $[ F ]$, however there is no
  other spherical class $Q$ with $[ Q ] \cdot [ Q ] = 0$ and such that $0 <
  \omega ( Q ) < \omega ( [ F ] )$. Thus $\gamma$ must fix $[ F ]$.

  Next we show that $i$ is a deformation retract.  Let $\psi\co D^n
  \rightarrow \sft$ such that $\psi ( \partial D^n ) \subset
  \mathcal{S}$. We will produce a retraction of $\psi$ to a disc of
  symplectomorphisms, while fixing its boundary.

  We do this by first producing a retraction of forms via inflation, and then
  applying Moser's Lemma:

  \begin{definition}
    Denote by $\mathcal{P}^{}$ the space of symplectic forms on $X$ which
    are cohomologous to $[ \omega ]$ and which restrict to symplectic forms,
    agreeing with the orientation induced by $\omega$, on each member of the
    triple $( F, Z_0, Z_{\infty} )$.
  \end{definition}

  Consider the sphere of symplectic forms $\psi^{\ast} ( \omega ) = \bigcup_{x
  \in D^n} \psi^{} ( x )^{\ast} ( \omega )$. Then $\psi^{\ast} ( \omega )
  \subset \mathcal{P}$. One can show by inflating along $Z_0$  that
  $\mathcal{P}$ is weakly contractible. The proof is similar to that in
  {\cite{Lalonde}} except that we must do it more parameters. We give a sketch
  below. Thus we can find a contraction of $\psi^{\ast} ( \sigma )$ to the
  constant sphere. Moser's Lemma then yields a family of diffeomorphisms
  \[ M_{\psi, t}\co D^n \times I \rightarrow \mathcal{D} ( M ) \]
  such that:
  \begin{enumerate}
    \item $M_{\psi, 1} ( d )^{\ast} ( \omega ) = \psi_d^{\star} ( \omega)$

    \item $M_{\psi, 0} ( d ) = \tmop{id}$

    \item $M_{\psi, t} ( \partial D ) = \tmop{id}$
  \end{enumerate}
  The map $M_{\psi, t}^{- 1} ( d ) \psi ( d )\co D^n \times I \rightarrow
  \sft$ then yields a retraction of $\psi$ into $\mathcal{S}$ as $t$
  travels from $0$ to $1$.

  Now we prove that $\mathcal{P}$ is weakly contractible.  Let $\phi\co
  S^n \rightarrow \mathcal{P}^{}$ be a sphere of symplectic forms based
  at $\omega$. Begin by homotoping $\phi$ so that it is constant in a
  neighborhood $U_b$ of the basepoint $b.$

  Since the set of closed 2--forms positive on each member of our triple is
  convex, $\phi$ admits a homotopy $\phi_t$ to the constant map through
  (possibly non-symplectic) closed 2--forms which are positive on the triple $(
  F_Z, Z_0, Z_{\infty} )$. Let $U$ be a neighborhood of $S^n \times 0$ such
  that each form $\sigma \in U$ is symplectic. Let
  \[ \chi\co S^n \times [ 0, 1 ] \rightarrow [ 0, 1 ] \]
  be a continuous function which vanishes on $U_0 \cup U_b \times [ 0, 1 ]$and
  is positive elsewhere.

  Denote by $\sigma_{\Sigma}$ an area form on $\Sigma$.  Then for some $\kappa
  \gg 0$ each form in $\phi_t + \kappa \chi ( x, t ) \pi_f^{\ast} (
  \sigma_{\Sigma} )$ is symplectic {\cite{thurstonsympfib}}, and the new
  homotopy
  \[ \begin{array}{ll}
       \phi^1_t ( x, t ) = \phi_t ( x, t ) +  \kappa \chi ( x, t
       ) \pi_f^{\ast} ( \sigma_{\Sigma} ) & x \in S^n, t \in [ 0, 1 ]
     \end{array} \]
  travels through symplectic forms, positive on our triple.

  However, we pay a price: the cohomology classes of the forms $
  \text{$\phi^1_t$} ( x, t )$ lie on the line $[ \omega ] + s \cdot
  \kappa [ F ]$, where $PD ( \cdot )$ denotes Poincare duality, and $s \in
  [0,1]$. We will now alter the homotopy $\phi^1_t$ by adding a sufficient
  multiples of a Thom class of the section $ Z_0$ for each value of $t$
  so that the forms in the new homotopy are Poincare dual to a multiple
  of  $[ \omega ] = aPD ( [ Z_0 ] ) + bPD ( [ F ] )$.  This process is
  called inflation. Its earliest appearance came in the papers of McDuff
  on the classification of symplectic structures on ruled surfaces. The
  version we will use is more refined:

  \begin{lemma}[McDuff {\cite{Sympalmostcomp}}]
    Let $M$ be a 4--manifold with a compact
    family of  of symplectic structures $\gamma$, and a compact family of
    almost complex structures $J_{\gamma}\co\Gamma\rightarrow\mathcal{J}$,
    which make a symplectic curve $C$ with $C \cdot C \geq 0$ holomorphic, and
    such that $\gamma$(x) tames $J_{\gamma} ( x )$ for each $x \in \Gamma$.

    Then for each $\beta > 0$ there is a compact family of closed 2--forms
    $\tau_{\gamma}\co\Gamma \rightarrow \Omega^2$, supported in an arbitrarily
    small neighborhood of $C$, and such that the form $\text{$\gamma$(x)} +
    \tau_{\gamma} ( x )$ is symplectic, tames $J_{\gamma} ( x )$ and has
    cohomology class $[ \omega ] + \beta PD ( [ C ] )$.
  \end{lemma}

  Applying Proposition \ref{lem Preserving two plane bundles}
  in the Appendix  we can find a {family $J_{\phi}\co S^n
  \times [ 0, 1 ] \rightarrow \mathcal{J}$ of almost complex structures on $M$
  such that $J_{\phi} ( x, t )$ is tamed by $\phi^1_t ( x )$, and which
  make each member of our triple holomorphic.} We then apply McDuff's Lemma
  above  $C = Z_0$ and $\beta = \frac{a \kappa}{b}$, and denote the resulting
  family of 2--forms by $\text{} \tau_{\phi} ( x, t )$. The set of forms
  taming a given almost complex structure is convex. Thus as $\phi^1_t ( x )$
  and $\text{} \phi^1_t ( x ) + \tau_{\phi} ( x, t )$ both tame
  $J_{\phi} ( x, t )$ so does $\phi^1_t ( x ) + s \cdot \tau_{\phi} ( x,
  t )$ for $s \in [ 0, 1 ]$. Thus all of these forms are positive on the
  triple $( F_Z, Z_0, Z_{\infty} )$.

  Consider, then, the homotopy $\phi_t^1$:
  $$\phi^2_t ( x ) = \phi^1_t ( x ) + \chi ( x, t ) \tau_{\phi}(x,t)$$
  Then
  \begin{eqnarray*}
    [ \phi_t^2 ( x ) ] & = & aPD ( [ Z_0 ] ) + bPD ( [ F_Z ] ) + \chi ( x, t )
    \kappa PD ( [ F_Z ] ) + \chi ( x, t ) \frac{a \kappa}{b} PD ( [ Z_0 ] )\\
    & = & ( 1 + {\kappa \chi ( x, t )}/{b} ) [ \omega ]
  \end{eqnarray*}
  Thus we have moved our homotopy to one which takes place only
  in classes which are multiples of $[ \omega ]$. One can then rescale each
  part of the homotopy by the appropriate constant factor to obtain a homotopy
  of our sphere within the original cohomology class
  $$[ \phi_t^3 ( x ) ] =
    ( \phi^2_t ( x ) + \chi ( x, t ) \tau_{\phi} ( x, t )/{( 1 + {\kappa \chi ( x, t )}/{b} )} \eqno{\qed}$$

\begin{remark}
  The proof of Proposition \ref{S in Soft a retract} does not require that
  $( X, \omega )$ satisfies the cohomology assumption.  However our later
  applications will require this assumption.
\end{remark}

\subsection{$ \sft$ ``acts transitively'' on $\mathcal{F Z}$}

We now show that the ``group'' $ \sft$ ``acts transitively'' on
$\mathcal{F Z} .^{}$  We shall see that, as a result, its homotopy
type is amenable to computation via $J$--holomorphic curves.

\begin{lemma}
  \label{lem:SF0inf =3D Z0infF}For every triple $( S_0, S_{\infty}, F_S ) \in
  \mathcal{F Z}$ there is a $g \in \mathcal{S}^{\tmop{soft}}$ such that $( S_0,
  S_{\infty}, F_S ) = ( g ( S_0 ), g ( S_{\infty} ), g ( F_S ) )$
\end{lemma}

\begin{proof}
  There is a diffeomorphism taking any symplectic fibrations by 2--spheres
  with fiber in class $[ F ]$ to any other (see
  {\cite{Sympalmostcomp}}). Call this diffeomorphism $g_F$.

  Next, note that as $Z_0$ and $S_0$ have the same self-intersection, their
  normal bundles are isomorphic. Thus one can find a fiber-preserving
  diffeomorphism $g_0$ of a neighborhood of $Z_0$ to $S_0$. Then, as the
  diffeomorphisms of a disc which fix its boundary form a contractible set, we
  can extend this to diffeomorphism-preserving $F_S$.

  $g_0 \circ g_F ( F, Z_0,^{} Z_{\infty} ) = ( F_S, S_0, g_0 \circ g_F (
  Z_{\infty} ) )$. We now aim to find a diffeomorphism $g_{\infty}$ which
  preserves $F_S$ and $S_0$ and carries $g_0 \circ g_F ( Z_{\infty} )$ to
  $S_{\infty}$. Denote the sections of $F_S$ which miss $S_0$ by
  $\mathcal{Z}_2$. These are sections of the disc bundle $F_S - S_0^{}$. Thus
  $\mathcal{Z}_2$ is contractible, and we may can find an isotopy of (possibly
  non-symplectic) sections from $g_0 \circ g_F ( Z_{\infty} )$ to $S_{\infty}$
  lying in $\mathcal{Z}_2$. This isotopy may then be induced by a path of
  diffeomorphisms which preserve both $F_S$ and $S_0$. Let $g_{\infty}$ be the
  end of this path of diffeomorphisms. $g_{\infty} \circ g_0 \circ g_F$ is
  then a diffeomorphism carrying $( F_Z, Z_0, Z_{\infty} )$ into $( F_S, S_0,
  S_{\infty} )$.
\end{proof}

\section{\label{Jhol on rational}$J$--holomorphic curves on rational surfaces}

In this subsection we supply the necessary background from the theory of
$J$--holomorphic curves on symplectic sphere bundles over surfaces. The main
geometric ingredient in our proof is the following Proposition:

\begin{proposition}[Gromov{\cite{phol}}, McDuff
  {\cite{Sympalmostcomp}}]
  \label{pro:(Gromov-Mcduff)}
  Let $( X, F, Z_0, \omega )$  satisfy the cohomology assumption (Definition
  \ref{def:cohomology assumption}). Then for every almost complex structure
  $J$ tamed by $\omega$ there is a $J$--holomorphic fibration by spheres in
  class $[ F ]$.
\end{proposition}

\begin{proof}
  This is well known for generic {\acs}s, and if the base $\Sigma$ is not a
  sphere. So suppose $\Sigma = S^2$, and let $J_n$ be a sequence of generic
  almost complex structures approximating $J$.

  We begin by showing that there is a section $Z$ of the fibration $F$ such
  that $[ Z ] \cdot [ Z ] = 0$ or 1, and $\omega [ Z ] \geq \omega [ F ]$.  Note
  that there are only two topological $S^2$--bundles over $S^2$.  One is
  trivial, and the other is not. For the trivial bundle, every section has
  even self-intersection, and for the non-trivial bundle every section
  has odd self-intersection.

  For the trivial bundle there is a section $Z$, of zero self-intersection.
  $[ Z ] = [ Z_0 ] - \frac{k}{2}$[F], and so $\omega [ Z ] = k( 1 -
  \varepsilon )/2$.  Thus, as $\omega ( F ) = ( 1 - \varepsilon )$ and  $k \geq
  2$, $\omega [ Z ] \geq \omega [ F ]$. For the non-trivial bundle there is a
  section $Z$ such that $[ Z ] \cdot [ Z ] = 1$.  Then $[ Z ] = [ Z_0 ] -
  \frac{k - 1}{2} [ F ]$, and $\omega ( Z ) = ( k + 1 ) ( 1 -
  \varepsilon )/2 \geq \omega ( F )$, since $k \geq 1$.  Let $x_1 \in F ( \Sigma
  )$.  Then for each $J_n$ there is a unique smooth irreducible curve $F_n$
  through $x_1$. I claim these cannot degenerate in the limit.

 For, suppose
  these degenerate to a cusp curve $B$.
  Denote the irreducible components of $B$ by $B_i$. Then $[ B_i ] = a_i [ Z
  ] + b_i [ F ]$, where $a_i$ and $b_i$ always have opposite signs.  Since
  $\sum_{i \in I} a_i = 0$, there must be one of each sign.  Denote the
  component with positive $a_i$ by $B_+$ and that with negative $a_i$ by
  $B_-$.

  Now choose a point $x_2 \in F ( \Sigma ) \backslash C_1$.  There cannot be
  an irreducible curve in class $[ F ]$ passing through $x_2$ as $[ F ] \cdot
  [ B_- ] < 0$.  Thus there must be a cusp curve $C$ passing through it.  It
  too splits into components $C_i$ with the above properties, with $C_+$ and
  $C_-$ defined analogously.  But then $[ C_+ ] \cdot [ B_+ ] < 0$,
  contradicting positivity of intersection.
\end{proof}

We will most often use it in the following form:

\begin{lemma}
  \label{Construction of F transv to pair}Let $( X, \omega )$
  be a symplectic 4--manifold  admitting a symplectic fibration by 2--spheres by
  $F$, which satisfies the cohomology assumption. Let $\{ S_i \}$be a
  collection of symplectic curves such that $[ S_i ] \cdot [ F ] = 1$, $J$ a
  tamed almost complex structure which preserves each curve. Then there is a
  $J$- holomorphic fibration by 2 spheres $F_J$ in the class of $[ F ]$, such
  that with respect to this fibration $F_J$ the curves $\{ S_i \}$ are
  symplectic sections.
\end{lemma}

\begin{proof}
  By Proposition \ref{pro:(Gromov-Mcduff)} there is a unique $J$--holomorphic
  fibration $F$ by 2 spheres in class $[ F ]$. As $J$ is tamed this fibration
  is symplectic. By positivity of intersection each fiber must meet each curve
  transversely, precisely once.
\end{proof}

\subsection{The main geometric lemma}

\begin{definition}
  Denote by $\mathcal{Z}$ the space of pairs $( S_0, S_{\infty} )$ of
  disjoint symplectic curves in $( X, \omega )$ such that $[ S_0 ] = [ Z_0 ]$
  and $[ S_{\infty} ] = [ Z_{\infty} ]$.
\end{definition}

Note that there is an inclusion $i\co \mathcal{E}_{\infty} \hookrightarrow \mathcal{Z}$,
given by identifying $\mathcal{E}_{\infty}$ with the pairs of disjoint
symplectic curves $( S_{0,} S_{\infty} )$ where $S_{\infty} = Z_{\infty}$.

\begin{proposition}[Main geometric lemma]
  \label{lem:Forgetting fibration is a fibration}
  The forgetful map $\pi\co \mathcal{F Z}^{} \rightarrow
  \mathcal{Z}$ is a fibration with contractible fiber.
\end{proposition}

\begin{proof}
  The proof of this proposition will rely heavily on the results in the
  Appendix on almost complex structures.

  We begin by showing that $\pi$ is a fibration. $\pi$ is surjective, for as
  the curves $Z_0$ and $Z_{\infty}$ are disjoint from one another we can find
  a tamed almost complex structure $J$ which makes each curve holomorphic.
  Lemma \ref{Construction of F transv to pair} then provides a fibration $F$
  so that $( F, Z_0, Z_{\infty} ) \in \mathcal{F Z}$.

  I claim that $\pi$ has path lifting: Let $B$ be a polyhedron. We consider
  $\Phi\co B \times I \rightarrow \mathcal{Z}$, along with a lifting
  $\Phi_\text{lift}\co B \times 0 \rightarrow S \mathcal{F}$. We aim to extend
  $\Phi_\text{lift}$ to all of $B \times I$. Applying Proposition
  \ref{lem Preserving two plane bundles} we can find a $\Phi^J\co B \times I
  \rightarrow \mathcal{J}$,  such that $\Phi ( b, t )$ is
  $\Phi^J(b,t)$--holomorphic, and $\Phi_\text{lift} ( b, 0 )$ is
  $\Phi^J(b,0)$--holomorphic.
  Applying Lemma \ref{Construction of F transv to pair} we gain a family of
  fibrations, $\Phi_\text{lift} ( b, t )$ extending our original lifting on $B
  \times 0$.

  Finally we show that $\pi$ has contractible fiber. Denote by $\mathcal{J}_Z$
  the tamed almost complex structures which make both $Z_0$ and $Z_{\infty}$
  holomorphic. It is enough to show that the map
  \begin{eqnarray*}
    \rho\co \mathcal{J}_Z & \rightarrow & \pi^{- 1} ( Z_0, Z_{\infty} )\\
    \rho ( J ) & = & ( F, Z_0, Z_{\infty} ),
  \end{eqnarray*}
  where $F$ is the unique $J$--holomorphic fibration determined by Lemma
  \ref{Construction of F transv to pair}, is a fibration with contractible
  fiber. For then $\rho$ will be a weak homotopy equivalence, and as
  $\mathcal{J}_Z$ is also contractible, so must $\pi^{- 1} ( Z_0, Z_{\infty}
  )$ be contractible. We commence with this task.

  We first show that $\rho$ is a fibration on its image, ie that it has
  path lifting: Let $B$ be a polyhedron. Consider $\Phi\co B \times
  I \rightarrow \pi^{- 1} ( Z_0, Z_{\infty} )$, along with a lifting
  $\Phi_\text{lift}\co B \times 0 \rightarrow \mathcal{J}_{0, \infty}$
  such that $\Phi ( b, 0 )$ is $\Phi_\text{lift} ( b, 0 )$--holomorphic. Then
  Proposition \ref{lem Preserving two plane bundles} allows us to extend
  $\Phi_\text{lift}$ to all of $B \times I$.

  Let $( F, Z_0, Z_{\infty} ) \in \pi^{- 1} ( Z_0, Z )$. Then $\rho^{- 1} ( F,
  Z_0, Z_{\infty} ) = \mathcal{J}_{\tmop{FZ}}$, the space of almost complex
  structures making each member of the triple holomorphic. However, again by
  applying Proposition \ref{lem Preserving two plane bundles},
  we see that $\mathcal{J}_{\tmop{FZ}}$ is non-empty and contractible. Thus
  $\rho$ is  surjective with contractible fiber.
\end{proof}

\section{Understanding the orbit of $\sft$: Proofs of Contractibility and
Transitivity} \label{sec:ContractibilityandTransitivity}

We now complete the proof of Theorem \ref{thm:Symp(M)is-homotopy-equivalent
LSigma} by proving Proposition \ref{rational surface}:

\begin{prop*} [Problem on rational surfaces]
The space
$\mathcal{S^{} ( \text{}} X )_{\overline{\infty}}$,
of symplectomorphisms
of $X$ which fix a framing of $Z_{\infty}$, is homotopy equivalent
to the space $\mathcal{E}_{\infty}$ of unparametrized, embedded
symplectic surfaces $S$ of $X$ disjoint from $Z_{\infty}$.
Moreover, when $\Sigma$ has genus $0$, $\mathcal{E}_{\infty}$ is
contractible.
\end{prop*}

We will show that $\mathcal{S^{} ( \text{}} X )_{\overline{\infty}}$  acts
transitively on $\mathcal{E}_{\infty}$, and that the stabilizer $\mathcal{S^{}
( \text{}} X )_{\overline{\infty}, 0}$ is contractible. Then the fibration:
\[ \text{$\mathcal{S^{} ( \text{}} X )_{\overline{\infty}, 0}$} \rightarrow
   \mathcal{S^{} ( \text{}} X )_{\overline{\infty}} \rightarrow
   \mathcal{E}_{\infty} \]
gives a homotopy equivalence between $\mathcal{S^{} ( \text{}} X
)_{\overline{\infty}}$ and $\mathcal{E}_{\infty}$.

\subsection{\label{sub:Proof-of-Proposition symp acts trans on miss
Zinft}$\mathcal{\mathcal{\mathcal{}}}
\mathcal{S}_{\overline{\infty^{}}}$ acts transitively on
$\mathcal{E}_{\infty}$}

\begin{proof}
  We begin by showing that $\mathcal{S}$ acts transitively on
  $\mathcal{Z}$. It is enough to show that there is a symplectomorphism
  carrying $( Z_0, Z_{\infty} )$ to any other pair $( Z_0^1, Z_{\infty}^1 )$
  in $\mathcal{Z}$. Let $J$ be an almost complex structure leaving both
  $Z_0$ and $Z_{\infty}$ invariant. Apply Lemma \ref{Construction of F transv
  to pair} and denote the resulting fibration by $F$. Then by Lemma
  \ref{lem:SF0inf =3D Z0infF} there is a $\alpha_1 \in \sft )$ which carries
  $( F, Z_0, Z_{\infty} )$ into $( F^1, Z_0^1, Z_{\infty}^1 )$. Since
  $\mathcal{S} \hookrightarrow \sft$ is a deformation retract by Proposition
  \ref{Symp a retract of comp Diff}, there is an isotopy $\alpha_t$ through
  $\sft$ to a symplectomorphism $\alpha_0$. Applying this isotopy to $( Z_0,
  Z_{\infty} )$ yields a path of pairs of curves $\alpha_t ( Z_0, Z_{\infty}
  )$ which begins at $\alpha_1 ( Z_0, Z_{\infty} ) = ( Z_0^1, Z_{\infty}^1 )$
  and ends at $\alpha_0 ( Z_0, Z_{\infty} )$ within the orbit of $( Z_0,
  Z_{\infty} )$ under $\mathcal{S} .$ One can then induce this path $\alpha_t
  ( Z_0, Z_{\infty} )$ by a path of symplectomorphisms $\Psi_t$, constructed
  by an easy application of Moser's Lemma. Then $\Psi_1 \alpha_0 ( Z_0,
  Z_{\infty} ) = ( Z_0^1, Z_{\infty}^1 )$.

  I claim that $\mathcal{E}_{\infty} \subset \mathcal{Z}$. $\mathcal{}
  \textrm{}$ $\mathcal{S}_{\overline{\infty}}$ are then precisely the
  symplectomorphisms which preserve $\mathcal{E}_{\infty}$ and act
  transitively on this space.

  We must show that if $Z$ is a symplectic curve in $\widehat{E}
  \backslash Z_{\infty}$, abstractly symplectomorphic to $Z_0$, then $[
  Z ] = [ Z_0 ]$. As $[ Z_0 ]$ and $[ F ]$ span $H_2 ( \widehat{E} )$,
  \[ [ Z ] = a [ Z_0 ] + b [ F ]. \]
  I claim that $b = 0$. For, as $Z$ misses $Z_{\infty}$,
  \begin{eqnarray*}
    0 & = & [ Z ] \cdot [ Z_{\infty} ]\\
    & = & ( a [ Z_0 ] + b [ F ] ) \cdot [ Z_{\infty} ]\\
    & = & a [ Z_0 ] \cdot [ Z_{\infty} ] + b [ F ] \cdot [ Z_{\infty} ]\\
    & = & 0 + b.
  \end{eqnarray*}
  Moreover, as the $Z$ and $Z_0$ are abstractly symplectomorphic,
  \begin{eqnarray*}
    \omega [ Z ] & = & \omega [ Z_0 ]\\
    a \omega [ Z_0 ] & = & \omega [ Z_0 ]
  \end{eqnarray*}
  $$a = 1,\leqno{\rm and\ thus}$$
  hence $[ Z ] = [ Z_0 ]$, and $\mathcal{\mathcal{\mathcal{E}}}_{\infty}
  \subset \mathcal{Z}$.
\end{proof}

\subsection{\label{sub:Proof-of-Proposition SYmp} $\mathcal{S^{} ( \text{}} X
)_{\overline{\infty}}$  is homotopy equivalent to  $\mathcal{E}_{\infty}$ }

Denote by $\mathcal{S}_{_{\infty},_0}$ the symplectomorphisms that
preserve both $Z_{\infty}$ and $Z_0$. Denote by $\left( \sft
\right)_{_{\infty},_0}$ the diffeomorphisms in $\sft$ which do the
same.

\begin{proposition}
  \label{pro:Symp fixing both homotopy equiv to Diff fixing both}
  $\sio \hookrightarrow \sftio$ is a homotopy equivalence.
\end{proposition}

\begin{proof}
  Since $\mathcal{S}$ acts transitively on $\mathcal{Z}$ (see
  Section~\ref{sub:Proof-of-Proposition symp acts trans on miss
Zinft}), the orbit map
  $\phi\co \mathcal{S} \rightarrow \mathcal{Z}$ is a fibration. Consider
  \[ \eta\co \mathcal{} \sft \rightarrow S \mathcal{F}^{} \rightarrow
     \mathcal{Z}; \]
 $\eta$ is the composition of two fibrations. The second, $S \mathcal{F}^{} \rightarrow
  \mathcal{Z}$, is a fibration by Proposition \ref{lem:Forgetting fibration
  is a fibration}; that the first is a fibration follows immediately from the definitions.   The fiber
  of $\eta$ is $\sftio$.

  The inclusion $\mathcal{S} \hookrightarrow \sft$ yields a morphism
  \[ \begin{array}{ccc}
       \sio & \overset{i_1}{\hookrightarrow} & \sftio\\
       \Big\downarrow &  & \Big\downarrow\\
       \mathcal{S} & \overset{i_2}{\hookrightarrow} & \sft\\
       \Big\downarrow \rlap{$\phi$} &  & \Big\downarrow \rlap{$\eta$}\\
       \mathcal{Z} & \overset{( \text{id} )}{\longrightarrow} & \mathcal{Z}
     \end{array} \]
  of fibrations:
  $i_2$ and $( \text{id} )$ are homotopy equivalences, thus so is $i_1$ by the
  5--Lemma (Lemma \ref{lem:5Lemma}).
\end{proof}

\subsubsection{We utilize $J$--holomorphic curves}

\begin{proposition}
  \label{pro:Diff same as those preserving fib} $\sftio$ is homotopy
  equivalent to $\mathcal{D}_{\infty, 0}$, the diffeomorphisms which preserve
  the fibration $F_Z$ and both sections $Z_0$ and $Z_{\infty}$
\end{proposition}

\begin{proof}
  Denote by $\mathcal{F} ( Z_0, Z_{\infty} )$ the subset of $S
  \mathcal{F}^{}$ given by triples $( F, S_0, S_{\infty} )$ where $S_0 =
  Z_0$ and $S_{\infty} = Z_{\infty}$. Restricting the fibration of $\sft
  \rightarrow \mathcal{F Z}^{}$ to $\mathcal{F} ( Z_0, Z_{\infty} )$ yields
  a fibration
  \[ \mathcal{D}_{\infty, 0} \rightarrow \sftio \rightarrow \mathcal{F (} Z_0,
     Z_{\infty} ). \]
  As $\mathcal{F (} Z_0, Z_{\infty} )$ is the fiber of the
  forgetful fibration $\pi\co S \mathcal{F}^{} \rightarrow \mathcal{Z}$
  it is contractible by the main geometric lemma (Proposition
  \ref{lem:Forgetting fibration is a fibration}).
\end{proof}

Denote the symplectomorphisms which preserve $Z_0$ and fix both $Z_{\infty}$
and its normal bundle by $\mathcal{S}_{\overline{\infty^{}}, 0}$.
Denote the fiber preserving diffeomorphisms which preserve $Z_0$, and fix both
$Z_{\infty}$ and its normal bundle by $\mathcal{D}_{0, \overline{\infty}}$
$\mathcal{\text{}}$.

\begin{proposition}
  \label{symp is homotopy equivalent to Diff}
  $\mathcal{S}_{\overline{\infty^{}}, 0}$ is homotopy equivalent to
  $\mathcal{D}_{\overline{\infty}, 0}$.
\end{proposition}

The proof is a fairly straightforward application of the previous ideas.  One
obtains a morphism of fibrations from the inclusions $\text{$\mathcal{}$}
\text{$\mathcal{\mathcal{\text{$\mathcal{D}_{0, \overline{\infty}}$ }}}$}
\hookrightarrow \sftio$, $\text{$\mathcal{S}_{\overline{\infty^{}}, 0}
\text{}$} \hookrightarrow \sftio$, and the ``actions'' of each.  (Those which are not
groups still have the analogous fibrations.)  Then we apply the 5--lemma and
Proposition \ref{pro:Diff same as those preserving fib}. For details see
{\cite{thesis}}.

\subsubsection{We use the Riemann mapping theorem in parameters}

The contractibility of $\text{$\mathcal{S}_{\overline{\infty^{}}, 0}$}$ now
follows from combining the above Proposition \ref{symp is homotopy equivalent
to Diff} with the following:

\begin{lemma}
  $\mathcal{D}_{ 0, \overline{\infty}}$ is
  contractible.
\end{lemma}

\begin{proof}
  As the elements in
  $\mathcal{D}_{0, \overline{\infty}}$ fix the section
  $Z_{\infty}$, and preserve the fibration, they must in fact preserve each
  fiber of $F_Z$. Thus $\mathcal{D}_{0, \overline{\infty}}$ is the space of sections of a bundle over $\Sigma$ whose
  fiber consists of the diffeomorphisms of $S^2$ which fix a point $0$ (where
  $Z_0$ intersects the fiber) and the neighborhood of another point $\infty$
  (where $Z_{\infty}$ intersects the fiber). This fiber is thus contractible
  by the Riemann mapping theorem. The space of sections$\mathcal{}$
  $\mathcal{}  \mathcal{D}_{0, \overline{\infty}}$ is thus
  also a contractible set.
\end{proof}

\subsection{If $\Sigma$ is a sphere,  $\mathcal{E}_{\infty} $ is contractible}

\begin{proposition}
  \label{Sphere is contractible}If $\Sigma$ is a sphere, $\mathcal{E}_{\infty}
  $ is contractible.
\end{proposition}

\begin{proof}
  Denote by $\mathcal{J}_{\infty}$ the set of tamed almost complex structures
  on $X$ which make$Z_{\infty}$ holomorphic. Denote by
  $\mathcal{J}_{\infty}^S$ the space of pairs$( J,S
  )$ where $J \in \mathcal{J}_{\infty}$ and $S \in
  \mathcal{E}_{\infty}$ such that $S$ is $J$--holomorphic. We remind the
  reader that $k = [ Z_0 ] \cdot [ Z_0 ] = [S] \cdot [S]$.

  I claim that the projections $\pi_{\Sigma}\co \mathcal{J}_{\infty}^S
  \rightarrow \mathcal{E_{\infty}}$ and $\pi_J\co \mathcal{J}_{\infty}^S
  \rightarrow \mathcal{J}_{\infty}$ are  fibrations with contractible fiber,
  and thus homotopy equivalences.  As $\mathcal{J}_{\infty}$ is contractible,
  this will show that $\mathcal{E}_{\infty}$ must also be contractible.

  We first prove the claim for $\pi_\Sigma$.  The fiber of $\pi_{\Sigma}$
  is the set of tamed complex structures which make both $Z_{\infty}$and $S$
  holomorphic. As $Z_{\infty}$and $S$ form a disjoint pair of symplectic
  curves this is a contractible set.

  We now prove the claim for $\pi_J$.  Denote the 2 disc by $D^2$,
  and let $J \in
  \mathcal{J}_{\infty}$. Fix $k + 1$ distinct points $x_i$ on $Z_{\infty}$. By
  Proposition \ref{pro:(Gromov-Mcduff)} there is a unique $J$--holomorphic
  curve $F_i$ in class $[ F ]$ which passes through $x_i$.As both $S$ and the
  $F_i$ are $J$--holomorphic they must intersect positively. Thus $S$ meets
  each $F_i$ in precisely one point $\sigma_i$. As $S$ misses $Z_{\infty}$,
  $\sigma_i \in F_i - x_i \simeq D^2$. Lemma \ref{lem:unique curve} below
  shows that for any $(k{+}1)$--tuple in $\prod_{i = 1\ldots k + 1}
  F_i - x_i$ there is a unique such curve $S$.  Therefore $\pi_J$ is a
  homotopy equivalence and $\mathcal{E}_{\infty} $ is contractible.
\end{proof}

\begin{lemma}
  \label{lem:unique curve}Let $J \in \mathcal{J}_{\infty}$, and let $F_i$
  denote $k + 1$ distinct holomorphic spheres in class $[ F ]$. If $\Sigma =
  S^2$, then for any $(k{+}1)$--tuple of points in $\prod_{i = 1\ldots
  k + 1} F_i - x_i$ there is a unique, smooth $J$--holomorphic sphere
  in $Z_0$ which passes through them.
\end{lemma}

\begin{proof}
  For a generic $J$ the moduli space of $J$--holomorphic spheres through $q$
  points has dimension
  \[  2 c_1 ( T ( \widehat{E} ) ) ( [ Z_0 ] ) - 2 q - 2 .\]
  $$2 c_1 ( T ( \widehat{E_{\epsilon}} ) ) ( [ Z_0 ] ) = 2 ( \chi ( Z_0 ) + [
     Z_0 ] \cdot [ Z_0 ] ) = 4 + 2 k\leqno{\rm Now}$$
  so for the dimension to be 0 we need
  \[ q = k + 1. \]
  The Gromov--Witten invariant for this class is 1. Thus there is a
  $J$--holomorphic curve $\Theta$ through any $k + 1$ points. I claim that this
  curve is unique. Let $\Theta_1$ and $\Theta_2$ be two curves through these $k
  + 1$ points. Then these two curves must coincide by positivity of
  intersection as $[ Z_0 ] \cdot [ Z_0 ] = k$.

  I claim that $\Theta$ is always smooth and irreducible. For one can always
  approximate $J$ by a sequence of complex structures $J_n$ so that the $J_n$
  holomorphic curve through these $q$ points $\Theta_n$ is smooth. The
  sequence of curves $\Theta_n$ then converges to $\Theta$, and $\Theta$ is
  thus controlled by Gromov compactness.  We will now show that the need to
  \begin{enumerate}
    \item intersect the curves in class $[ F ]$ positively (curves in class
    $[ F ]$ exist for every $J$ tamed by $\omega$ by Proposition
    \ref{pro:(Gromov-Mcduff)}), and

    \item intersect $Z_{\infty}$ positively ($J \in J_{\infty}$ and thus
    $Z_{\infty}$ is a $J_{\infty}$ holomorphic curve)
  \end{enumerate}
  eliminate all such nodal curves, save those of the form
  \[ Z_{\infty} \cup \bigcup_{i = 1}^k F_i \]
  where the $F_i$ are (possibly repeated) spheres in class $F$. However curves
  of this last form are eliminated as well. They have only $k$ fiber curves
  $F_i$, they cannot pass through all $k + 1$ points. For  each point lies off
  $Z_{\infty}$ and in a distinct $J$--holomorphic fiber.  Note that this
  argument fails when the genus of $\Sigma > 0$.  For then, $q$, the number of
  points required by the dimension formula, is less than  $k$ and there is no
  contradiction.

  \paragraph{All nodal curves are of the form $Z_{\infty} \cup \bigcup_{i =
  1}^k F_i$}

  The second homology of $\widehat{E}_{\epsilon}$ is spanned by $[ Z_0 ]$
  and the fiber class $[ F ]$. The class of each irreducible component
  $\Theta_i$ of a curve may thus be written $a_i [ Z_0 ] + b_i [ F ]$. Each
  $a_i > 0$ as $a_i = [ \Theta_i ] \cdot [ F ]$ and each of the classes is
  represented by a holomorphic curve.

  The union of these components lies in class $[ Z_0 ]$, and thus
  \[ \sum_i ( a_i [ Z_0 ] + b_i [ F ] ) = [ Z_0 ]. \]
  As all the $a_i$ are positive integers, the only possibility which remains
  is that one $a_i = 1$ and the rest vanish. Moreover, for all $i$ such that
  $a_i = 0$, $b_i$ must be positive, as $\omega$ evaluated on each component
  must be positive. Thus we have reduced ourselves to:
  \begin{enumerate}
    \item A curve in class $[ Z_0 ] - l [ F ]$ for $l \in \mathbb{Z}, l > 0.$

    \item A collection of curves in class $b_i [ F ]$ $b_i > 0$ such that
    $\sum_i b_i = l$.
  \end{enumerate}
  Since there is a unique curve through each point in class $[ F ]$ these
  curves of ``type 2'' must be unions of fibers in $F$. Positivity of
  intersection with $Z_{\infty}$ then implies that the only $J$--holomorphic
  curve in class $[ Z_0 ] - l [ F ]$ ($l > 0$) is $Z_{\infty}$ itself with $l
  = k$.
\end{proof}

\section{Appendix}

\subsection{Tamed almost complex structures preserving sub-bundles}

In this subsection we collect the results we require about tamed almost
complex structures preserving sub-bundles. They are listed below in order of
their difficulty. The first two are classical, the last less so, and we
provide a proof of it here.

\begin{definition}
  If $\pi\co V \rightarrow B$ is a symplectic vector bundle, let $\pi_J\co
  \mathcal{J} ( V ) \rightarrow B$ be the bundle such that $\pi_J^{- 1} ( b )$
  are the tamed almost complex structures on $\pi^{- 1} ( b )$. If $\eta_i
  \subset V$ are symplectic sub bundles, let $\pi_J\co\mathcal{J} ( V, \eta_1,
  \eta_2, \ldots ) \rightarrow B$ be the (possibly locally non-trivial) bundle
  where $\pi_J^{- 1} ( b )$ are the tamed almost complex structures on $\pi^{-
  1} ( b )$, which preserve each $\eta_i .$

\end{definition}

The first result goes back at least to Gromov's seminal paper {\cite{phol}}.

\begin{lemma}
  Let $( V, \omega ) \rightarrow B$ be a symplectic vector bundle over a
  polyhedron $B$. Then $\rho\co\mathcal{J} ( V ) \rightarrow B$ is a bundle
  with contractible fibers.
\end{lemma}

Next we consider almost complex structures preserving a given plane field.
This is also a classical result:

\begin{lemma}
  \label{lem:preservin one plane bundle contr}Let $( V, \omega ) \rightarrow
  B$ be a 4--dimensional symplectic vector bundle over a polyhedron $B$. Let
  $\vartheta \subset V$ be a 2--dimensional symplectic sub-bundle of $V$. Then
  $\rho\co\mathcal{J} ( V, \vartheta ) \rightarrow B$ is a bundle with
  contractible fibers.
\end{lemma}

Let $Q \subset B$ be a sub-polyhedron of $B$. Then given a section $\phi_Q$ of
$\rho\co\mathcal{J} ( V, \vartheta ) \rightarrow Q$, we may construct
a section $\phi$ of $\rho\co\mathcal{J} ( V, \vartheta ) \rightarrow B$
extending $\phi_Q$.

Finally we will need to consider the tamed almost complex structures
preserving two transverse plane fields. Preserving two planes requires a good
deal more work than preserving one, as pairs of symplectic planes have moduli.
This is less well known and we include a proof of it here.

\begin{proposition}
  \label{lem Preserving two plane bundles}Let $( V, \omega ) \rightarrow B$ be
  a 4--dimensional symplectic vector bundle over a polyhedron $B$. Let
  $\vartheta_1, \vartheta_2 \subset V$ be 2--dimensional symplectic sub-bundles
  of $V$, such that $\vartheta_1, \vartheta_2$ intersect transversely in each
  fiber, and the the symplectic orthogonal projection $\pi_{12}^{\perp}\co
  \vartheta_1 \rightarrow \vartheta_2$ is orientation-preserving. Let $Q
  \subset B$ be a sub-polyhedron, and let $\phi_Q$ be a section of $\rho\co
  \mathcal{J} ( V, \vartheta_1, \vartheta_2 ) \rightarrow B$, defined over
  $Q$.
  Then there is a section $\phi$ of $\rho$ which extends $\phi_Q$.
\end{proposition}

\begin{proof}
  Constructing $\phi$ is equivalent to constructing a section of $\phi^1
  \oplus \phi^2$ of $J ( \vartheta_1 ) \oplus J ( \vartheta_2 )$ such that the
  resulting almost complex structure is tamed by $\omega$. Denote by
  $\phi_Q^i$ the sections such that
  \[ \phi_Q = \phi_Q^1 \oplus \phi_Q^2. \]
  Constructing $\phi_1$ alone is fairly simple, for by
Lemma~\ref{lem:preservin one plane bundle contr},
  \[ J ( \vartheta_1 ) \rightarrow B \]
  is a bundle with contractible fibers. Thus this bundle admits a section
  $\phi^1$ extending $\phi_Q^1$. We now proceed with the problem of
  constructing $\phi^2$ extending $\phi_Q^2$.

  Our main tool in will be the following lemma in linear algebra, which
  provides the almost complex structures satisfying our conditions with a
  convex structure. This will allow to use partitions of unity to construct
  $\phi^2$.

  \begin{lemma}
    \label{lem:Linear algebra}Let $V, P$ be two symplectic planes in $R^4$ with
    symplectic structure $\omega$. Let $\pi_{\perp}\co P \rightarrow V^{\perp}$
    denote the symplectic orthogonal projection. Suppose that $\pi_{\perp}$ is
    orientation-preserving. Fix an $\omega -$tame complex structure $J_P$ on
    $P$, and $v \in V$, with $v \neq 0$. Denote the space of almost complex
    structures $J_V$ on $V$ such that $J = ( J_P \oplus J_V )$ is
    $\omega$--tame by $\mathcal{J}_V$.

    Then $\Phi_v\co \mathcal{J}_V \rightarrow V$, given by $J_V \rightarrow J_V
    ( v )$, gives a homeomorphism of $\mathcal{J}_V$ onto a convex set.
  \end{lemma}

  The proof of this lemma is involved, and we defer it to the end of this
  subsection. For now we concentrate on its application to our argument. It
  has the following immediate consequence:

  \begin{lemma}
    Let $s$ be a section of $\vartheta_2$, non-vanishing over a set $X_s
    \subset B$. Then there is a section $\phi^2 \in \mathcal{J} ( \vartheta_2
    )$ such that $\phi^2$ extends $\phi_Q^2$ and $\phi^1 \oplus \phi^2$ is
    tamed by $\omega$ over $X_s$.
  \end{lemma}

  \begin{proof}
    As $s$ is non-vanishing, $\alpha$ is determined by its action on $s$.
    Lemma \ref{lem:Linear algebra} tells us that the set of allowable choices
    for $\alpha ( s )$ form an open, convex set. As $X_s$ is paracompact, so
    we can use a partition of unity to construct a section $s_{\alpha}$ of
    $\vartheta_2$ over $X_s$, so that for
    \begin{align*}
      \phi^2&\co s \longrightarrow s_{\alpha}\\
      \phi^2&\co s_{\alpha} \longrightarrow - s
    \end{align*}
    the almost complex structure $\phi^1 \oplus \phi^2$ is tamed by $\omega$.
  \end{proof}

  Proposition \ref{lem Preserving two plane bundles} then follows by applying
  a partition of unity to a covering $X_{s_i}$ coming from a finite set of
  sections $\{ s_i \}$ of $\Phi_Z^{1 \ast} ( \eta )$ such that $\bigcup_i
  X_{s_i} = \Sigma \times P$.
\end{proof}

\subsubsection{Proof of the linear algebra lemma}

In this subsection we provide the proof of the promised linear
algebra lemma \ref{lem:Linear algebra}

\begin{lem*}[Linear algebra lemma]  Let $V, P$ be two
symplectic planes in $R^4$ with symplectic structure $\omega$. Let
$\pi_{\perp}\co P \rightarrow V^{\perp}$ denote the symplectic
orthogonal projection. Suppose that $\pi_{\perp}$ is
orientation-preserving. Fix an $\omega -$tame complex structure $J_P$
on $P$.
Denote the space of almost complex structures $J_V$ on $V$ such
that $J = ( J_P \oplus J_V )$ is $\omega$--tame by $\mathcal{J}_V$.

Then for each non-zero $v \in V$, $\Phi_v\co \mathcal{J}_V
\rightarrow V$ given by $J_V \rightarrow J_V ( v )$ gives a
homeomorphism of $\mathcal{J}_V$ onto a convex set.
\end{lem*}

\begin{proof}

 We begin by establishing some useful coordinates.   Let $\pi\co P \rightarrow V^{\perp}$ denote the
symplectic
  orthogonal projection to $V$.  We choose a $v_1 \in \tmop{im} ( \pi^{} )$,
  and denote $\pi^{- 1} ( v_1 )$ by $p$.  Write
  \[ \label{eq:definition of p} p = w_1 + v_1 \]
  where $w_1 \in V^{\perp}$. Write
  \[ \label{eq:defintion of jp} J_P ( p ) = w_2 + v_2 \]
  where $w_2 \in V^{\perp}, v_2 \in V$. Then
  \[ J ( w_1 ) = w_2 + v_2 - J_V ( v_1 ). \]
  Applying $J$ to both sides of this equation we compute $Jw_2$:
  \[ J ( w_2 ) = - w_1 - J_V ( s_2 ) - s_1 \]
  Throughout this lemma we will suppress $\omega$ and just denote the pairing
  of two vectors $p$ and $q$ by $( p, q ) .$

  $\pi$ is orientation-preserving. Thus, as $( w_1, J ( w_2 ) ) > 0$,  so
  is $( w_1, w_2 ) > 0$. To lessen our burden of constants, scale $v_1$, thus
  scaling $p = \pi^{- 1} ( v_1 )$, so that
  \[ ( w_1, w_2 ) = 1 \]
  This scaling in turn dilates the image of $\Phi_{v_1}$, and thus does not
  affect its convexity.

  As $P$ is symplectic $( w_1, w_2 ) + ( v_1, v_2 ) > 0$ and thus
  \begin{equation}
    ( v_1, v_2 ) > - 1. \label{eq:bound on (s1,s2)}
  \end{equation}
  \paragraph{Reduction to Cauchy--Schwartz}
  We now commence in earnest. Let $w + v \in V^{\perp} \oplus V = R^4$. What must we require of $J_V$ so
that $( w + v, J ( w + v ) ) > 0$ for all such pairs $w$ and $v$?
  \begin{eqnarray*}
    ( w + v, J ( w + v ) ) & = & ( w, Jw ) + ( w, Jv ) + ( v, Jw ) + ( v, Jv )
  \end{eqnarray*}
  $( w, Jv ) = 0$ as $J$ must preserve $V$. And we have
  \[ ( w, Jw ) + ( v, Jv ) + ( v, Jw ). \]
  The first two terms are positive. Also $( v, Jw ) = ( v, q )$ where $q$ is the
  projection of $Jw$ to $V$; this term may well be negative. We seek to bound
  its absolute value in terms of the other two (positive) terms.

  We replace $v$ by $- Jv$ throughout the equation. As we seek a bound for
  all pairs $w, v$ this has no effect on our task. Moreover as $( v, Jv ) = (
  - JJv, - Jv )$ this has no effect on the third term. The second term $( v, q
  )$ becomes $( - J_V v, q ) = ( v, J_V q ) .$

  We seek to show that
  \[ \left| ( v, J_V q ) \right| \leq ( w, Jw ) + ( v, Jv ). \]
  As $J_V$ has determinate 1, it preserves $\omega |_V$, thus $( \cdot, J_V
  \cdot )$ is a (symmetric) inner product on $V$. Cauchy--Schwartz then implies
  that
  \[ | ( v, J_V q )| \leq ( v, J_V q )^{\frac{1}{2}} ( q, J_V q
     )^{\frac{1}{2}}. \]
  If $v = kq$, this bound is achieved. Thus the tamed $J_V$ are
  precisely those such that
  \begin{equation}
    ( v, J_V v )^{\frac{1}{2}} ( q, J_V q \dot{)}^{\frac{1}{2}} < ( w, Jw ) +
    ( v, Jv ). \label{eq:inequality from CS}
  \end{equation}
  \paragraph{Understanding the constraint imposed by Cauchy--Schwartz}\nl
  We now unpack
  this inequality. Write $w$ as $aw_1 + bw_2$. Then
  \begin{eqnarray*}
    ( w, Jw ) & = & ( aw_1 + bw_2, J ( aw_1 + bw_2 ) )\\
    & = & \left( aw_1 + bw_2, a ( w_2 + v_2 - J_V v_1 ) + b ( - w_1 - J_V v_2
    - v_1 ) \right)\\
    & = & a^2 ( w_1, w_2 ) - b^2 ( w_2, w_1 )\\
    & = & a^2 + b^2
  \end{eqnarray*}
and
  $$( q, J_V q ) = ( av_2 - aJ_V v_1 - bJ_V v_2 - bv_1, aJ_V s_2 + av_1 +
    bv_2 - bJ_V v_1 ).$$
  Expanding the right hand side creates sixteen pairings, however some of them are
  zero, and the four ``$ab$'' terms all cancel. Upon summing we are left with
  \begin{eqnarray*}
    ( q, J_V q ) & = & \lambda ( v_1, J_V s_1 ) + \lambda ( v_2, J_V v_2 ) - 2
    \lambda ( v_1, v_2 )
  \end{eqnarray*}
  where we denote $( a^2 + b^2 )$ by $\lambda$.

  Our inequality \eqref{eq:inequality from CS} then reads
  $$( v, J_V v )^{\frac{1}{2}} ( \lambda ( v_1, J_V v_1 ) + \lambda (v_2,J_V
    v_2 ) - 2 \lambda ( v_1, v_2 ) )^{\frac{1}{2}} < \lambda + (v,J_V v).$$
  If $v = 0$ the inequality places no restriction on $J_V .$ Thus we may
  assume that $v$ is not zero. Since the condition $( w + v, J ( w + v ) ) >
  0$ is invariant under scaling by a positive constant, we may scale the
  vector $v + w$ so that $( v, J_V v ) = 1$, if we assume that $J_V$ tames
  $\omega$ on $V.$ We do so, and are left with one free parameter $\lambda>0$.
  $$( \lambda ( v_1, J_V v_1 ) + \lambda ( v_2, J_V v_2 ) - 2 \lambda ( v_1,
    v_2 ) )^{\frac{1}{2}} < \lambda + 1.$$
  Squaring both sides yields
  $$\lambda | ( v_1, J_V v_1 ) + ( v_2, J_V v_2 ) - 2 ( v_1, v_2 ) | <
    \lambda^2 + 2 \lambda + 1.$$
  This may be achieved for all $\lambda$ if and only if
  \begin{equation}
    | ( v_1, J_V v_1 ) + ( v_2, J_V v_2 ) - 2 ( v_1, v_2 ) | < 4
    \label{eq:4bound}
  \end{equation}
  At this point our proof bifurcates into three cases.

 \textbf{Case 1: $\tmop{rk} ( \pi ) = 2$}\qua
  We may assume that $v = v_1$, and
  write $J_V v_1 = cv_1 + dv_3$. We now describe the constraints that
  \eqref{eq:4bound} places on $c$ and $d$.
  \[ J_V v_2 = - \frac{c^2+1}{d} v_1 - cv_2 \]
  Substituting into \eqref{eq:4bound} we get
  \[ | ( d + (c^2 + 1)/d - 2 ( v_1, v_2 ) | < 4. \]
  $J_V$'s tameness restricted to $V$ \footnote{We assumed this when we scaled
  $v$ so that $( v, Jv ) = 1$.} translates to $d$ having the same sign as $(
  s_1, s_2 )$. Thus $( d + (c^2 + 1)/d )$ and $- 2 ( s_1, s_2 )$ have
  opposite sign, and our inequality is equivalent to
  \[ | ( d + (c^2 + 1)/d ) | < 4 + 2 ( v_1, v_2 ). \]
  If we denote $4 + 2 ( s_1, s_2 )$ by $\gamma$, the set of solutions of
  this inequality form a disc centered at
  \[ ( c, d ) = ( 0, \gamma/2 ) \]
  with radius $( \gamma^2/4 - 1 )^{\frac{1}{2}}$. Since $( v_1, v_2 )
  > - 1$ by \eqref{eq:bound on (s1,s2)}, $\gamma > 2$ and this disc is
  non-empty.

  \textbf{Case 2: $\tmop{rk} ( \pi ) = 1$}\qua
  In this case $v_1$ and $v_2$ are
  linearly dependent. So choose a $v_3 \in V$ which is independent of $V_1$,
  such that $( v_1, v_3 ) = 1$. Moreover if $v \neq v_1$ we choose $v_3$ so
  that $v = k v_3$.
  Then
  \begin{align*}
  J_V v_1 & = cv_1 + dv_3 \\
  J_V v_3 & = - \frac{c^2+1}{d} v_1 - cv_3.
  \end{align*}
  As $v_2 = k v_1$, the constraints that \eqref{eq:4bound} places on $c$ and $d$
  are much weaker. Substituting into \eqref{eq:4bound} we get
    $$|c| < \frac{4}{1 + k^2}$$
  and no condition on $d$. Furthermore, $v$ is either $v_1$ or $v_3$, and in both cases the
  resulting sets of vectors $J_V v_1$ (when $v = v_1$) and $J_V v_3$ (when $v
  = v_3$) form the band $cv_1 + dv_3$ where
    $$|c| < \frac{4}{1 + k^2}.$$
 Thus they form a convex set.

\textbf{Case 3: $\tmop{rk} ( \pi ) = 0$ }\qua
  This case actually requires no proof at all, for here the two planes are orthogonal. Thus for any $J_V$
tamed by $\omega |_V$, $J_V \oplus J_P$ is tamed by $\omega$, and
$\tmop{im} ( \Phi_v )$ is
  a half plane.
\end{proof}

\subsection{A non-traditional 5--Lemma}

The 5--Lemma is usually presented in the context of chain
complexes, and as such it is usually stated as a Lemma about
Abelian groups. However its usual proof actually applies in much
more generality. As we will require it, we present the more general
statement here. The proof is the standard one which the reader can
find in {\cite{Spanier}}.

\begin{lemma}
  \label{lem:5Lemma}
  Let
  \[ \begin{CD}
  G_5 @>\alpha_5>> G_4 @>\alpha_4>> G_3 @>\alpha_3>> G_2 @>\alpha_2>> G_1 \\
  @V{\gamma_5}VV @V{\gamma_4}VV @V{\gamma_3}VV @V{\gamma_2}VV @V{\gamma_1}VV \\
  H_5 @>\beta_5>> H_4 @>\beta_4>> H_3 @>\beta_3>> H_2 @>\beta_2>> H_1
  \end{CD} \]
  be a diagram of pointed sets, with each row exact. Suppose that $G_3$ and
  $G_2$ are groups, $\gamma_i$ makes $H_i$ a $G_i$--set, and the morphisms
  $\alpha_3$ and $\beta_3$ respect this structure. Suppose further that
  $\gamma_1$, $\gamma_2$, $\gamma_4$ and $\gamma_5$ are bijections.  Then
  $\gamma_3$ is a bijection.
\end{lemma}

We will most often use the above generalized 5--lemma in the
following form, gained by applying it to the long exact sequence
of homotopy groups of a fibration. One needs the extra generality
to deal with morphisms on $\pi_0$ and $\pi_1$

\begin{lemma}
  \label{lem:FibrationLem}
  Let
  \[ \begin{CD}
  G_0 @>\phi_0>> F_i \\
  @VVV           @VVV \\
  G_1 @>\phi_1>> F_1 \\
  @VVV           @VVV \\
  G_2 @>\phi_2>> F_2 \\
  \end{CD} \]
  be a morphism of fibrations such that:
  \begin{enumerate}
    \item $G_i$ is a group.

    \item $H_i$ is a $G_i$--set.

    \item The maps $\phi_1$ and $\phi_2$ are  weak homotopy equivalences.
  \end{enumerate}
  Then $\phi_0$ is also a homotopy equivalence.
\end{lemma}

\end{document}